\documentclass{amsart}%
\usepackage{amsfonts}
\usepackage{amsmath}
\usepackage{amssymb}
\usepackage{graphicx}%
\newtheorem{theorem}{Theorem}
\theoremstyle{plain}

\newtheorem{corollary}{Corollary}

\newtheorem{definition}{Definition}
\newtheorem{example}{Example}

\newtheorem{lemma}{Lemma}

\newtheorem{proposition}{Proposition}

\numberwithin{equation}{section}
\begin{document}
\title[Singular polynomials]{Singular polynomials and modules for the symmetric groups}
\author{Charles F. Dunkl}
\address{Department of Mathematics, University of Virginia\\
P.O.Box 400137 Charlottesville VA 22904-4137}
\email{cfd5z@virginia.edu}
\urladdr{http://www.people.virginia.edu/\symbol{126}cfd5z/}
\thanks{During the preparation of this paper the author was partially supported by NSF
grant DMS 0100539.}
\date{January 26, 2005}

\begin{abstract}
For certain negative rational numbers $\kappa_{0}$, called singular values,
and associated with the symmetric group $S_{N}$ on $N$ objects, there exist
homogeneous polynomials annihilated by each Dunkl operator when the parameter
$\kappa=\kappa_{0}$. It was shown by de Jeu, Opdam and the author (Trans.
Amer. Math. Soc. 346 (1994), 237-256) \ that the singular values are exactly
the values $-\frac{m}{n}$ with $2\leq n\leq N$, $m=1,2,\ldots$ and $\frac
{m}{n}$ is not an integer. For each pair $\left(  m,n\right)  $ satisfying
these conditions there is a unique irreducible $S_{N}$-module of singular
polynomials for the singular value $-\frac{m}{n}$. The existence of these
polynomials was previously established by the author (IMRN 2004, \#67,
3607-3635). The uniqueness is proven in the present paper. By using Murphy's
(J. Alg. 69(1981), 287-297) results on the eigenvalues of the Murphy elements,
the problem of existence of singular polynomials is first restricted to the
isotype $\tau$ (where $\tau$ is a partition of $N$ corresponding to an
irreducible representation of $S_{N}$) satisfying the condition that
$n/\gcd\left(  m,n\right)  $ divides $\tau_{i}+1$ for $1\leq i<l$; $l$ is the
length of $\tau$, that is, $\tau_{l}>\tau_{l+1}=0$. Then by arguments
involving the analysis of nonsymmetric Jack polynomials it is shown that the
assumption $\tau_{2}\geq n/\gcd\left(  m,n\right)  $ leads to a contradiction.
This shows that the singular polynomials are exactly those already determined,
and are of isotype $\tau$, where $\tau_{2}=\ldots=\tau_{l-1}=\left(
n/\gcd\left(  m,n\right)  \right)  -1\geq\tau_{l}$.

\end{abstract}
\maketitle

\section{Introduction}

The symmetric group $S_{N}$ on $N$ letters acts on $\mathbb{R}^{N}$ by
permutation of coordinates. The alternating polynomial, also called the
discriminant, is defined by $a_{N}\left(  x\right)  =\prod\nolimits_{1\leq
i<j\leq N}\left(  x_{i}-x_{j}\right)  $ for $x\in\mathbb{R}^{N}$ and is a
fundamental object associated to the group. The Macdonald-Mehta-Selberg
integral for $S_{N}$ is%
\[
\left(  2\pi\right)  ^{-N/2}\int_{\mathbb{R}^{N}}\left\vert a_{N}\left(
x\right)  \right\vert ^{2\kappa}\exp\left(  -\frac{1}{2}\sum_{i=1}^{N}%
x_{i}^{2}\right)  dx=\prod\limits_{n=2}^{N}\frac{\Gamma\left(  n\kappa
+1\right)  }{\Gamma\left(  \kappa+1\right)  },
\]
for $\kappa\geq0$. The right hand side is a meromorphic function of $\kappa$
without zeroes and with poles at $\kappa=-\frac{m}{n}$, for $2\leq n\leq
N,m=1,2,3,\ldots$ and $\frac{m}{n}$ is not an integer. (For an algebraic proof
of the integral, see \cite[Sect. 8.7]{DX}.) Do these values have another
connection with the symmetric group? The purpose of this paper is to show that
for each pair $\left(  m,n\right)  $ of natural numbers with $2\leq n\leq N$
and $\frac{m}{n}$ not an integer there is a unique irreducible $S_{N}$-module
of homogeneous polynomials which have a certain singularity property with
respect to a commutative algebra of differential-difference operators. In a
previous paper \cite{D2} the author established the existence of a space of
such polynomials for each pair $\left(  m,n\right)  $. This paper proves the
uniqueness of the polynomials and the associated modules. By use of the Murphy
elements one can find a link between the singular polynomials, the partition
of $N$ which labels the module and the nonsymmetric Jack polynomials (NSJP's).
This is the family of simultaneous eigenvectors of a commuting set $\left\{
\mathcal{U}_{i}\left(  \kappa\right)  :1\leq i\leq N\right\}  $ of operators
(involving a parameter $\kappa$). The singular polynomials come from the
specializations of certain of NSJP's when $\kappa$ takes the value $-\frac
{m}{n}$. The algebra generated by the $\mathcal{U}_{i}\left(  \kappa\right)  $
is semisimple (that is, the set of NSJP's forms a basis for all polynomials)
for generic $\kappa$, but this property may be lost for some negative rational
values. A part of the development is to show how to find limits of certain
expressions in the NSJP's as $\kappa$ approaches $-\frac{m}{n}$.

Murphy \cite{Mu} found the eigenvalues of the Murphy elements when restricted
to any irreducible $S_{N}$-module. In Section 2 we use his results to find a
necessary condition on a partition to allow corresponding singular polynomials
and also to prove a uniqueness result. The condition is this: suppose
$\gcd\left(  m,n\right)  =1$ and there is an $S_{N}$-module of singular
polynomials corresponding to $\kappa=-\frac{m}{n}$, and suppose the module is
labeled by the partition $\tau$ (that is, $\tau=\left(  \tau_{1},\tau
_{2},\ldots\right)  $ with $\sum_{i\geq1}\tau_{i}=N$ and $\tau_{1}\geq\tau
_{2}\geq\ldots\geq0$) then $n|\left(  \tau_{i}+1\right)  $ for $1\leq
i<\ell\left(  \tau\right)  $, where $\ell\left(  \tau\right)  =\max\left\{
j:\tau_{j}\geq1\right\}  $). Section 3 develops the relevant results on
NSJP's. In Section 4 it is shown that the two-part partitions of the form
$\left(  N-n,n\right)  $ with $n|\left(  N-n+1\right)  $ can not give rise to
singular polynomials. Section 5 completes the proof of the main result: if
$\tau_{2}\geq n$ then there can not be corresponding singular polynomials.
This shows that the class of partitions, namely, $\tau$ with $n|\left(
\tau_{1}+1\right)  $ and $\tau_{2}=\tau_{3}=\ldots=\tau_{\ell\left(
\tau\right)  -1}=n-1$, appearing in \cite{D2} is exhaustive. In Section 6
there is a restatement of the main theorem and a discussion of the relation
between singular polynomials and modules over the rational Cherednik algebra.

The group $S_{N}$ is the finite reflection group of type $A_{N-1}$ and it
acts\ by permutation of coordinates. Let $\mathbb{N}_{0}$ denote $\left\{
0,1,2,3,\ldots\right\}  $ (also $\mathbb{N}:\mathbb{=}\left\{  1,2,3,\ldots
\right\}  $, $\mathbb{Z}$ and $\mathbb{Q}$ denote the sets of integers and
rational numbers respectively). For $\alpha\in\mathbb{N}_{0}^{N}$ (called a
\textquotedblleft composition\textquotedblright) let $\left\vert
\alpha\right\vert =\sum_{i=1}^{N}\alpha_{i}$ and define the monomial
$x^{\alpha}$ to be $\prod_{i=1}^{N}x_{i}^{\alpha_{i}}$; its degree is
$\left\vert \alpha\right\vert $. The length of a composition $\alpha$ is
$\ell\left(  \alpha\right)  =\max\left\{  j:\alpha_{j}>0\right\}  $. Consider
elements of $S_{N}$ as permutations on $\{1,2,\ldots,N\}$. Then, for
$x\in\mathbb{R}^{N}$ and $w\in S_{N}$ let $\left(  xw\right)  _{i}=x_{w\left(
i\right)  }$ for $1\leq i\leq N$ and extend this action to polynomials by
$\left(  wf\right)  \left(  x\right)  =f\left(  xw\right)  $. This has the
effect that monomials transform to monomials: $w\left(  x^{\alpha}\right)
=x^{w\alpha}$ where $\left(  w\alpha\right)  _{i}=\alpha_{w^{-1}\left(
i\right)  }$ for $\alpha\in\mathbb{N}_{0}^{N}$. (Consider $x$ as a row vector,
$\alpha$ as a column vector, and $w$ as a permutation matrix, with $1$'s at
the $\left(  w\left(  j\right)  ,j\right)  $ entries.) The reflections in
$S_{N}$ are the transpositions interchanging $x_{i}$ and $x_{j}$ and are
denoted by $\left(  i,j\right)  $ for $i\neq j$.

In \cite{D1} the author constructed for each finite reflection group a
parametrized commutative algebra of differential-difference operators. Let
$\kappa$ be a formal parameter, that is, $\mathbb{Q}\left(  \kappa\right)  $
is a transcendental extension of $\mathbb{Q}$.

\begin{definition}
The space of polynomials is $\mathcal{P}:=\mathrm{span}_{\mathbb{Q}\left(
\kappa\right)  }\left\{  x^{\alpha}:\alpha\in\mathbb{N}_{0}^{N}\right\}  $ and
for $n\in\mathbb{N}_{0}$ the subspace of homogeneous polynomials of degree $n$
is $\mathcal{P}_{n}:=$\newline$\mathrm{span}_{\mathbb{Q}\left(  \kappa\right)
}\left\{  x^{\alpha}:\alpha\in\mathbb{N}_{0}^{N},\left\vert \alpha\right\vert
=n\right\}  $. For $p\in\mathcal{P}$ and $\alpha\in\mathbb{N}_{0}^{N}$ let
$\mathrm{coef}\left(  p,\alpha\right)  $ denote the coefficient of $x^{\alpha
}$ in $p$ (thus $p=\sum_{\beta}\mathrm{coef}\left(  p,\beta\right)  x^{\beta}$).
\end{definition}

For the symmetric group $S_{N}$ the operators are defined as follows:

\begin{definition}
For any polynomial $f$ on $\mathbb{R}^{N}$ and $1\leq i\leq N$ let
\[
\mathcal{D}_{i}\left(  \kappa\right)  f\left(  x\right)  =\frac{\partial
}{\partial x_{i}}f\left(  x\right)  +\kappa\sum_{j\neq i}\frac{f\left(
x\right)  -\left(  ij\right)  f\left(  x\right)  }{x_{i}-x_{j}}.
\]

\end{definition}

It was shown in \cite{D1} that $\mathcal{D}_{i}\left(  \kappa\right)
\mathcal{D}_{j}\left(  \kappa\right)  =\mathcal{D}_{j}\left(  \kappa\right)
\mathcal{D}_{i}\left(  \kappa\right)  $ for $1\leq i,j\leq N$ and each
$\mathcal{D}_{i}\left(  \kappa\right)  $ maps $\mathcal{P}_{n}$ to
$\mathcal{P}_{n-1}$ for $n\geq1$. A specific numerical parameter value
$\kappa_{0}$ is said to be a \textit{singular value} (associated with $S_{N}$)
if there exists a nonzero polynomial $p$ such that $\mathcal{D}_{i}\left(
\kappa_{0}\right)  p=0$ for $1\leq i\leq N$ : such a $p$ is called a
\textit{singular polynomial}. It was shown in \cite{DJO} that the singular
values are the numbers $-\frac{m}{n}$ where $n=2,\ldots,N,\,m\in\mathbb{N}$
and $\frac{m}{n}\notin\mathbb{Z}$. Earlier, Opdam \cite{O} showed that the
$S_{N}$-Bessel function $J\left(  x,y\right)  $ considered as a function of
the parameter $\kappa$ has poles precisely at these numbers (for $\kappa>0$
the Bessel function is the entire solution of the system of equations
$\sum\limits_{j=1}^{N}\left(  \mathcal{D}_{j}^{\left(  x\right)  }\left(
\kappa\right)  \right)  ^{k}J\left(  x,y\right)  =\left(  \sum\limits_{j=1}%
^{N}y_{j}^{k}\right)  J\left(  x,y\right)  ,~1\leq k\leq N,~J\left(
0,y\right)  =1,~J\left(  xw,y\right)  =J\left(  x,yw\right)  =J\left(
x,y\right)  $ for $x,y\in\mathbb{C}^{N}$). Because the operators
$\mathcal{D}_{i}\left(  \kappa\right)  $ preserve homogeneity and have the
$S_{N}$-transformation property $\mathcal{D}_{i}\left(  \kappa\right)  \left(
i,j\right)  =\left(  i,j\right)  \mathcal{D}_{j}\left(  \kappa\right)  $, the
set of singular polynomials for a specific singular value is a direct sum of
irreducible $S_{N}$-modules of homogeneous polynomials. The set of partitions
of length $\leq N$ is denoted by $\mathbb{N}_{0}^{N,P}$ and consists of all
$\lambda\in\mathbb{N}_{0}^{N}$ such that $\lambda_{i}\geq\lambda_{i+1}$ for
$1\leq i\leq N-1$. When writing partitions it is customary to suppress
trailing zeros and to use exponents to indicate multiplicity, for example
$\left(  5,2^{3}\right)  $ is the same as $\left(  5,2,2,2,0\right)
\in\mathbb{N}_{0}^{5,P}.$ (The exponent notation is also used for compositions.)

The irreducible representations of $S_{N}$ are labeled by partitions $\tau$ of
$N$ (that is, $\tau\in\mathbb{N}_{0}^{N,P}$ and $\left\vert \tau\right\vert
=N$) and we say a polynomial $f$ is of \textit{isotype} $\tau$ if $f$ is an
element of an irreducible $S_{N}$-submodule of $\mathcal{P}_{n}$ corresponding
to $\tau$. It was conjectured in \cite{DJO} that the two-part representations
$\left(  n-1,N-n+1\right)  $ (with $2\left(  n-1\right)  \geq N$) give rise to
singular polynomials for the singular values $-\frac{m}{n}$ with $\gcd\left(
m,n\right)  <\frac{n}{N-n+1}$, and the representations $\left(
dn-1,n-1,\ldots,n-1,\tau_{l}\right)  $ for $d,n\in\mathbb{N}$ give rise to
singular polynomials for the singular values $-\frac{m}{n}$ with $\gcd\left(
m,n\right)  =1$ (where $l=\ell\left(  \tau\right)  $ and $N=\left(
dn-1\right)  +\left(  l-2\right)  \left(  n-1\right)  +\tau_{l}$). This
construction is presented in \cite{D2} in terms of nonsymmetric Jack
polynomials. In this paper we show that there are no other singular
polynomials. By using Murphy's techniques in his construction of the Young
seminormal representations \cite{Mu} we can show that the isotype $\tau$ of
any irreducible module of singular polynomials for $\kappa_{0}=-\frac{m}{n}$
(with $\gcd\left(  m,n\right)  =1$) must satisfy $n|\left(  \tau_{i}+1\right)
$ for $1\leq i<\ell\left(  \tau\right)  $. After that most of the work is to
show that the assumption $\tau_{2}\geq n$ leads to a contradiction. Note that
the condition $\tau_{2}<n$ implies for three or more parts $\tau_{i}=n-1$ for
$2\leq i<\ell\left(  \tau\right)  $ (and $\tau_{\ell\left(  \tau\right)  }\leq
n-1$) and for two parts that $\tau_{2}<\dfrac{\tau_{1}+1}{\gcd\left(
m_{1},\tau_{1}+1\right)  }=n$ for the singular value $-\frac{m_{1}}{\tau
_{1}+1}$. These are the restrictions described above.

The notation is almost the same as that in \cite{D2} except that the parameter
has been incorporated. Key parts of the proofs depend on the behavior of
polynomials as $\kappa$ approaches a singular value $\kappa_{0}$. The
commutative algebra of the operators defining the nonsymmetric Jack
polynomials is generated by
\[
\mathcal{U}_{i}\left(  \kappa\right)  f\left(  x\right)  =\mathcal{D}%
_{i}\left(  \kappa\right)  x_{i}f\left(  x\right)  -\kappa\sum_{j=1}%
^{i-1}\left(  j,i\right)  f\left(  x\right)  ,1\leq i\leq N.
\]
(this differs by the additive constant $\kappa$ from the notation in
\cite[Ch.8]{DX}). The operators act in a triangular manner on monomials, as is
explained below.

\begin{definition}
\label{order}For $\alpha\in\mathbb{N}_{0}^{N}$, let $\alpha^{+}$ denote the
unique partition such that $\alpha^{+}=w\alpha$ for some $w\in S_{N}$. For
$\alpha,\beta\in\mathbb{N}_{0}^{N}$ the partial order $\alpha\succ\beta$
($\alpha$ dominates $\beta$) means that $\alpha\neq\beta$ and $\sum_{i=1}%
^{j}\alpha_{i}\geq\sum_{i=1}^{j}\beta_{i}$ for $1\leq j\leq N$; $\alpha
\vartriangleright\beta$ means that $\left\vert \alpha\right\vert =\left\vert
\beta\right\vert $ and either $\alpha^{+}\succ\beta^{+}$ or $\alpha^{+}%
=\beta^{+}$ and $\alpha\succ\beta$. The notations $\alpha\succeq\beta$ and
$\alpha\trianglerighteq\beta$ include the case $\alpha=\beta$.
\end{definition}

When acting on the monomial basis of $\mathcal{P}_{n}$ the operators
$\mathcal{U}_{i}\left(  \kappa\right)  $ have on-diagonal coefficients
involving the following \textit{rank}\ function on $\mathbb{N}_{0}^{N}$. We
denote the cardinality of a set $E$ by $\#E$.

\begin{definition}
\label{rankdef}For $\alpha\in\mathbb{N}_{0}^{N}$ and $1\leq i\leq N$ let
\begin{align*}
r\left(  \alpha,i\right)   &  =\#\left\{  j:\alpha_{j}>\alpha_{i}\right\}
+\#\left\{  j:1\leq j\leq i,\alpha_{j}=\alpha_{i}\right\}  ,\\
\xi_{i}\left(  \alpha;\kappa\right)   &  =\left(  N-r\left(  \alpha,i\right)
\right)  \kappa+\alpha_{i}+1.
\end{align*}

\end{definition}

Clearly for a fixed $\alpha\in\mathbb{N}_{0}^{N}$ the values $\left\{
r\left(  \alpha,i\right)  :1\leq i\leq N\right\}  $ consist of all of
$\left\{  1,\ldots,N\right\}  $, are independent of trailing zeros (that is,
if $\alpha^{\prime}\in\mathbb{N}_{0}^{M},\alpha_{i}^{\prime}=\alpha_{i}$ for
$1\leq i\leq N$ and $\alpha_{i}^{\prime}=0$ for $N<i\leq M$ then $r\left(
\alpha,i\right)  =r\left(  \alpha^{\prime},i\right)  $ for $1\leq i\leq N$),
and $\alpha\in\mathbb{N}_{0}^{N,P}$ if and only if $r\left(  \alpha,i\right)
=i$ for all $i$ (the latter property motivated the use of \textquotedblleft%
$1\leq j\leq i$\textquotedblright\ rather than \textquotedblleft$1\leq
j<i$\textquotedblright\ in the definition). Then (see \cite[p.291]{DX})
$\mathcal{U}_{i}\left(  \kappa\right)  x^{\alpha}=\xi_{i}\left(  \alpha
;\kappa\right)  x^{\alpha}+q_{\alpha,i}\left(  x\right)  $ where $q_{\alpha
,i}\left(  x\right)  $ is a sum of terms $\pm\kappa x^{\beta}$ with
$\alpha\vartriangleright\beta$. The nonsymmetric Jack polynomials are the
simultaneous eigenvectors of $\left\{  \mathcal{U}_{i}\left(  \kappa\right)
:1\leq i\leq N\right\}  $ and they are well-defined for generic $\kappa$.

\section{$S_{N}$-modules}

In this section we find necessary conditions on a partition $\tau$ of $N$ for
the existence of singular polynomials of isotype $\tau$. Suppose $f$ is a
singular polynomial for some singular value $\kappa_{0}$. We may assume $f$ is
homogeneous because the operators $\mathcal{D}_{i}\left(  \kappa\right)  $ are
homogeneous and that $f$ has rational coefficients ($\mathcal{D}_{i}\left(
\kappa_{0}\right)  $ is a rational operator). Any translate of $f$ by $S_{N}$
is singular so $\mathrm{span}_{\mathbb{Q}}\left\{  wf:w\in S_{N}\right\}  $ is
an $S_{N}$-module of singular polynomials for $\kappa_{0}$. Suppose one of the
irreducible components has isotype $\tau$, for some partition $\tau\ $with
$\left\vert \tau\right\vert =N$. This decomposition is a computation over
$\mathbb{Q}$ (from the representation theory of $S_{N}$). Henceforth we
restrict our attention to this module, denoted by $M$.

We turn to the application of Murphy's results. For any given isotype he
determined the eigenvalues and transformation properties of the eigenvectors
of the commuting operators $\left\{  \sum\nolimits_{j=1}^{i-1}\left(
i,j\right)  :2\leq i\leq N\right\}  $ (Jucys-Murphy elements). The results
have to be read in reverse in a certain sense.

\begin{proposition}
Suppose $f$ is a singular polynomial for $\kappa=\kappa_{0}\in\mathbb{Q}$ and
$1\leq i\leq N$, then $\mathcal{U}_{i}\left(  \kappa_{0}\right)
f=f+\kappa_{0}\sum_{j=i+1}^{N}\left(  i,j\right)  f$.
\end{proposition}

\begin{proof}
We have the commutation $\mathcal{D}_{i}\left(  \kappa\right)  \left(
x_{i}f\right)  =x_{i}\mathcal{D}_{i}\left(  \kappa\right)  f+f+\kappa
\sum_{j\neq i}\left(  i,j\right)  f$. Now set $\kappa=\kappa_{0}$ and note
that $\mathcal{U}_{i}\left(  \kappa_{0}\right)  f=\mathcal{D}_{i}\left(
\kappa_{0}\right)  \left(  x_{i}f\right)  -\kappa_{0}\sum_{j<i}\left(
i,j\right)  f=-\kappa_{0}\sum_{j<i}\left(  i,j\right)  f$.
\end{proof}

Denote the Murphy elements $\omega_{i}=\sum\limits_{j=N-i+2}^{N}\left(
N+1-i,j\right)  $ for $2\leq i\leq N$ and let $\omega_{1}=0$ (as a
transformation); then $\mathcal{U}_{i}\left(  \kappa_{0}\right)
f=f+\kappa_{0}\omega_{N+1-i}f$ for $f\in M$. A \textit{standard Young tableau}
(SYT) of shape $\tau$ is a one-to-one assignment of the numbers $\left\{
1,\ldots,N\right\}  $ to the nodes of the Ferrers diagram $\left\{  \left(
i,j\right)  \in\mathbb{N}^{2}:1\leq i\leq\ell\left(  \tau\right)  ,1\leq
j\leq\tau_{i}\right\}  $ so that the entries increase in each row and in each
column. The notation $T\left(  i,j\right)  $ refers to the entry at row $i$,
column $j$. There is an order on SYT's of given shape (for details see
\cite[p.288]{Mu}) and the maximum SYT in this order, denoted by $T_{0}$, is
produced by entering the numbers $1,2,\ldots,N$ row by row (the first row is
$1,\ldots,\tau_{1}$, the second is $\tau_{1}+1,\ldots,\tau_{1}+\tau_{2}$ and
so forth).

\begin{definition}
Let $\tau\in\mathbb{N}_{0}^{N,P}$ with $\left\vert \tau\right\vert =N$, and
let $Y\left(  \tau\right)  $ denote the set of SYT's of shape $\tau$. Suppose
$T\in Y\left(  \tau\right)  $ then let $rw\left(  i,T\right)  ,cm\left(
i,T\right)  ,\eta_{i}\left(  T\right)  :=cm\left(  i,T\right)  -rw\left(
i,T\right)  $ denote the row, column and content, respectively of the node of
$T$ containing $i$, for $1\leq i\leq N$.
\end{definition}

(With this notation $T\left(  rw\left(  i,T\right)  ,cm\left(  i,T\right)
\right)  =i$.) Murphy constructed a basis $\left\{  f_{T}:T\in Y\left(
\tau\right)  \right\}  $ for the irreducible representation of isotype $\tau$
such that $\omega_{i}f_{T}=\eta_{i}\left(  T\right)  f_{T}$ for each $i$ and
$T$ (actually, this is an isomorphic image of the construction, which is in
terms of specific polynomials, of minimal degree). The eigenvalues $\left(
\eta_{1}\left(  T\right)  ,\ldots,\eta_{N}\left(  T\right)  \right)  $
determine the SYT $T$ uniquely thus there is a unique (up to scalar
multiplication) basis $\left\{  f_{T}:T\in Y\left(  \tau\right)  \right\}  $
for $M$ with%
\[
\mathcal{U}_{i}\left(  \kappa_{0}\right)  f_{T}=\left(  1+\kappa_{0}%
\eta_{N+1-i}\left(  T\right)  \right)  f_{T},\text{ for }1\leq i\leq N,~T\in
Y\left(  \tau\right)  .
\]
(The argument for uniqueness of $T$ is in \cite{Mu}: one can reconstruct $T$
by adjoining boxes containing $2,3,\ldots,N$ to $1$ by using the values
$\eta_{2}\left(  T\right)  ,\ldots,\eta_{N}\left(  T\right)  $; at any stage
the locations at which one can adjoin a box to make a larger SYT have
different contents.) We will show that $f_{T_{0}}$ is (a multiple of)
$x^{\lambda}+\sum_{\beta\vartriangleleft\lambda}A_{\beta}x^{\beta}$ with
coefficients $A_{\beta}\in\mathbb{Q}$ and $\lambda\in\mathbb{N}_{0}^{N,P}$
with%
\begin{equation}
\lambda_{N+1-i}=-\kappa_{0}\sum_{j=1}^{s-1}\left(  \tau_{j}+1\right)  ,\text{
for }\sum_{j=1}^{s-1}\tau_{j}<i\leq\sum_{j=1}^{s}\tau_{j}.\label{tau2lb}%
\end{equation}
This implies that if $\kappa_{0}=-\frac{m}{n}$ with $\gcd\left(  m,n\right)
=1$ then $n|\left(  \tau_{j}+1\right)  $ for $1\leq j<\ell\left(  \tau\right)
$ (the maximum value for $s$ in the above formula). The proof relies on the
triangularity properties of the $\mathcal{U}_{i}\left(  \kappa_{0}\right)  $
with respect to the order $\vartriangleright$.

\begin{definition}
For each $T\in Y\left(  \tau\right)  $ let $C_{T}=\left\{  \beta\in
\mathbb{N}_{0}^{N}:\mathrm{coef}\left(  f_{T},\beta\right)  \neq0\right\}  $.
Let $C$ be the set of $\alpha\in\cup_{T\in Y\left(  \tau\right)  }C_{T}$ such
that $\alpha$ is $\vartriangleright$-maximal in some $C_{T}$ (that is,
$\alpha,\beta\in C_{T}$ and $\beta\trianglerighteq\alpha$ implies
$\beta=\alpha$).
\end{definition}

\begin{lemma}
If $\alpha$ is a $\vartriangleright$-maximal element of $C$ then $\alpha$ is a partition.
\end{lemma}

\begin{proof}
It suffices to show that for any $\beta\in\cup_{T}C_{T}$ there exists a
partition $\lambda\in C$ with $\lambda\trianglerighteq\beta$. Since
$M=\mathrm{span}_{\mathbb{Q}}\left\{  f_{T}:T\in Y\left(  \tau\right)
\right\}  $ is $S_{N}$-invariant we see that for any $w\in S_{N}$ and
$\beta\in C_{T}$ for some $T$ there exists $T_{1}\in Y\left(  \tau\right)  $
such that $w\beta\in C_{T_{1}}$ (note that $wf_{T}\left(  x\right)
=f_{T}\left(  xw\right)  $ and $w\left(  x^{\beta}\right)  =x^{w\beta}$). So
$\cup_{T}C_{T}$ is $S_{N}$-invariant, in particular if $\beta\in\cup_{T}C_{T}$
then $\beta^{+}\in\cup_{T}C_{T}$. Thus there exists $\gamma\in C$ such that
$\gamma\trianglerighteq\beta^{+}$. Since $C$ is finite the maximal elements
$\alpha$ satisfy $\alpha\trianglerighteq\alpha^{+}$, that is, $\alpha$ is a partition.
\end{proof}

The next step is to show that there is a unique maximal element in $C$
determined by equation \ref{tau2lb}.

\begin{lemma}
\label{stdtab}If $T\in Y\left(  \tau\right)  $ satisfies $\eta_{s+1}\left(
T\right)  \leq\eta_{s}\left(  T\right)  +1$ for $1\leq s<N$ then $T=T_{0}.$
\end{lemma}

\begin{proof}
We have to show that the condition implies $rw\left(  s,T\right)  \leq
rw\left(  s+1,T\right)  $ for each $s$. Fix $s$ and let $T\left(  i_{1}%
,j_{1}\right)  =s$ and $T\left(  i_{2},j_{2}\right)  =s+1$ so that $\eta
_{s}\left(  T\right)  -\eta_{s+1}\left(  T\right)  =\left(  j_{1}%
-j_{2}\right)  +\left(  i_{2}-i_{1}\right)  $. We list the possibilities for
these nodes in any SYT. If $s$ and $s+1$ are in the same row of $T$ then
$i_{2}=i_{1},j_{2}=j_{1}+1$ and $\eta_{s+1}\left(  T\right)  =\eta_{s}\left(
T\right)  +1$. If $s$ and $s+1$ are in the same column of $T$ then
$i_{2}=i_{1}+1,j_{2}=j_{1}$ and $\eta_{s+1}\left(  T\right)  =\eta_{s}\left(
T\right)  -1$. The condition $i_{1}<i_{2}$ and $j_{1}<j_{2}$ is impossible or
else $s<T\left(  i_{2},j_{1}\right)  <s+1$. Also the condition $i_{1}>i_{2}$
and $j_{1}>j_{2}$ is impossible or else $s+1<T\left(  i_{1},j_{2}\right)  <s$.
If $i_{1}<i_{2}$ and $j_{1}>j_{2}$ then $\eta_{s}\left(  T\right)  -\eta
_{s+1}\left(  T\right)  \geq2$. The case $i_{1}>i_{2}$ and $j_{1}<j_{2}$ is
ruled out by hypothesis because it implies $\eta_{s}\left(  T\right)
-\eta_{s+1}\left(  T\right)  \leq-2$.
\end{proof}

\begin{theorem}
Suppose $\lambda$ is a $\vartriangleright$-maximal element of $C$ then
$\lambda$ is $\vartriangleright$-maximal in $C_{T_{0}}$ and is given by
equation \ref{tau2lb}.
\end{theorem}

\begin{proof}
By hypothesis $\lambda$ is a partition and is $\vartriangleright$-maximal in
$C_{T}$ for some $T\in Y\left(  \tau\right)  $. By the triangularity property
of $\mathcal{U}_{i}\left(  \kappa_{0}\right)  $ we have that
\begin{align*}
\mathrm{coef}\left(  \left(  1+\kappa_{0}\eta_{N+1-i}\left(  T\right)
\right)  f_{T},\lambda\right)   &  =\mathrm{coef}\left(  \mathcal{U}%
_{i}\left(  \kappa_{0}\right)  f_{T},\lambda\right) \\
&  =\xi_{i}\left(  \lambda;\kappa_{0}\right)  \mathrm{coef}\left(
f_{T},\lambda\right)  ,
\end{align*}
and $\xi_{i}\left(  \lambda;\kappa_{0}\right)  =\left(  N-i\right)  \kappa
_{0}+\lambda_{i}+1$, for $1\leq i\leq N$. This gives the equations%
\begin{align*}
\left(  N-i\right)  \kappa_{0}+\lambda_{i}+1  &  =1+\kappa_{0}\eta
_{N+1-i}\left(  T\right)  ,\\
\lambda_{N+1-i}  &  =\kappa_{0}\left(  \eta_{i}\left(  T\right)  +1-i\right)
.
\end{align*}
Since $\lambda$ is a partition $\lambda_{N+1-i}\leq\lambda_{N-i}$ for $1\leq
i<N$ and thus $\eta_{i}\left(  T\right)  +1-i\geq\eta_{i+1}\left(  T\right)
+1-\left(  i+1\right)  $ (note that $\kappa_{0}<0$). By Lemma \ref{stdtab}
$T=T_{0}$. By definition of $T_{0}$ for $1\leq i\leq\tau_{1}$ we have
$\eta_{i}\left(  T_{0}\right)  =i-1$ thus $\lambda_{N+1-i}=0$. In the range
$\sum_{j=1}^{s-1}\tau_{j}+1\leq i\leq\sum_{j=1}^{s}\tau_{j}$ (row $s$ of
$T_{0})$ $\eta_{i}\left(  T_{0}\right)  =\left(  i-\sum_{j=1}^{s-1}\tau
_{j}\right)  -s$ and $\lambda_{N+1-i}=-\kappa_{0}\left(  \sum_{j=1}^{s-1}%
\tau_{j}+s-1\right)  =-\kappa_{0}\sum_{j=1}^{s-1}\left(  \tau_{j}+1\right)  $.
\end{proof}

\begin{corollary}
\label{bigtau}There is a unique $\vartriangleright$-maximal element $\lambda$
of $C$ given by equation \ref{tau2lb}, and $n|\left(  \tau_{j}+1\right)  $ for
$1\leq j<\ell\left(  \tau\right)  $ (where $\kappa_{0}=-\frac{m}{n}$ and
$\gcd\left(  m,n\right)  =1$).
\end{corollary}

\begin{proof}
The uniqueness is now obvious. The equation $\lambda_{N+1-i}=m\sum
\limits_{j=1}^{s-1}\dfrac{\tau_{j}+1}{n}$ for $\sum\limits_{j=1}^{s-1}\tau
_{j}<i\leq\sum\limits_{j=1}^{s}\tau_{j}$ shows inductively that $n|\left(
\tau_{j}+1\right)  $ for $1\leq j<\ell\left(  \tau\right)  $; since the
maximum value of $s$ is $\ell\left(  \tau\right)  $.
\end{proof}

\begin{corollary}
For any $\kappa_{0}=-\frac{m}{n}$ and partition $\tau$ of $N$ there is at most
one irreducible $S_{N}$-module, consisting of singular polynomials for the
singular value $\kappa_{0}$, that has isotype $\tau$.
\end{corollary}

\begin{proof}
Suppose there are two unequal modules $M$ and $M^{\prime}$ satisfying the
hypotheses. Let $\left\{  f_{T}:T\in Y\left(  \tau\right)  \right\}  $ and
$\left\{  f_{T}^{\prime}:T\in Y\left(  \tau\right)  \right\}  $ be the
respective bases for $M$ and $M^{\prime}$ produced by Murphy's construction.
Normalize the two bases so that both $f_{T_{0}}$ and $f_{T_{0}}^{\prime}$ are
monic in $x^{\lambda}$ (that is, $f_{T_{0}}=x^{\lambda}+\sum_{\beta
\vartriangleleft\lambda}A_{\beta}x^{\beta}$ and $f_{T_{0}}^{\prime}$ has the
same form with $A_{\beta}$ replaced by $A_{\beta}^{\prime}$), with $\lambda$
given by equation \ref{tau2lb}. Let $g_{T}=f_{T}-f_{T}^{\prime}$ for $T\in
Y\left(  \tau\right)  $, then $\mathrm{span}_{\mathbb{Q}}\left\{  g_{T}:T\in
Y\left(  \tau\right)  \right\}  $ consists of singular polynomials and its
basis has the same transformation properties under the action of $S_{N}$ as
the basis of $M$. By the Theorem $\mathrm{coef}\left(  g_{T_{0}}%
,\lambda\right)  \neq0$, which is a contradiction.
\end{proof}

The following summarizes the results of this section. The polynomial
$f_{T_{0}}$ is renamed $g_{\lambda}$.

\begin{theorem}
\label{glambda}Suppose there exist singular polynomials for $\kappa_{0}%
=-\frac{m}{n}$ with $\gcd\left(  m,n\right)  =1$ (and $2\leq n\leq N$) of
isotype $\tau$, a partition of $N$, then $n|\left(  \tau_{i}+1\right)  $ for
$1\leq i<\ell\left(  \tau\right)  $ and there is a unique singular polynomial
$g_{\lambda}=x^{\lambda}+\sum_{\beta\vartriangleleft\lambda}A_{\beta}x^{\beta
}$ (with $A_{\beta}\in\mathbb{Q}$) of isotype $\tau$ such that $\mathcal{U}%
_{i}\left(  \kappa_{0}\right)  g_{\lambda}=\xi_{i}\left(  \lambda;\kappa
_{0}\right)  g_{\lambda}$ for $1\leq i\leq N$, where $\lambda$ is given by
equation \ref{tau2lb}.
\end{theorem}

\section{Nonsymmetric Jack polynomials}

These polynomials are the simultaneous eigenvectors of the commuting set of
operators $\left\{  \mathcal{U}_{i}\left(  \kappa\right)  :1\leq i\leq
N\right\}  $ . The existence follows from the triangular property and the fact
that the correspondence (from compositions to eigenvalues) $\alpha
\longmapsto\left(  \xi_{i}\left(  \alpha;\kappa\right)  \right)  _{i=1}^{N}$
\ is one-to-one for generic $\kappa$. We use the notation from \cite{D2} (for
now just the $x$-monic version is used but there will be a reference to the
$p$-monic version).

\begin{definition}
For $\alpha\in\mathbb{N}_{0}^{N}$, let $\zeta_{\alpha}^{x}\left(
\kappa\right)  $ denote the $x$-monic simultaneous eigenvectors, that is,
$\mathcal{U}_{i}\left(  \kappa\right)  \zeta_{\alpha}^{x}\left(
\kappa\right)  =\xi_{i}\left(  \alpha;\kappa\right)  \zeta_{\alpha}^{x}\left(
\kappa\right)  $ for $1\leq i\leq N$ and $\zeta_{\alpha}^{x}\left(
\kappa\right)  =x^{\alpha}+\sum\limits_{\beta\vartriangleleft\alpha}%
A_{\beta\alpha}^{x}\left(  \kappa\right)  x^{\beta}$, with coefficients
$A_{\beta\alpha}^{x}\left(  \kappa\right)  \in\mathbb{Q}\left(  \kappa\right)
$.
\end{definition}

\begin{definition}
\label{defbij}For $1\leq i\leq N$ the operators $\mathcal{B}_{ij}$ (with
$j\neq i$) and the operator $\mathcal{B}_{i}$ (each maps $\mathcal{P}_{n}$
into itself, for $n\in\mathbb{N}_{0}$) are given by%
\begin{align*}
\mathcal{B}_{ij}p\left(  x\right)   &  :=\frac{x_{i}p\left(  x\right)
-x_{j}p\left(  x\left(  i,j\right)  \right)  }{x_{i}-x_{j}}-\left\{
\begin{array}
[c]{ll}%
0, & i<j\\
p\left(  x\left(  i,j\right)  \right)  , & i>j,
\end{array}
\right. \\
\mathcal{B}_{i}p  &  :=\sum_{j\neq i}\mathcal{B}_{ij}p,\text{ for }%
p\in\mathcal{P}.
\end{align*}

\end{definition}

In this notation $\mathcal{U}_{i}\left(  \kappa\right)  p\left(  x\right)
=\frac{\partial}{\partial x_{i}}\left(  x_{i}p\left(  x\right)  \right)
+\kappa\mathcal{B}_{i}p\left(  x\right)  $. There is an easily proved
identity: $\mathcal{B}_{ij}+\mathcal{B}_{ji}=1$, and this shows directly that%
\[
\sum\limits_{i=1}^{N}\mathcal{U}_{i}\left(  \kappa\right)  =N+\sum
\limits_{i=1}^{N}x_{i}\dfrac{\partial}{\partial x_{i}}+\kappa\frac{N\left(
N-1\right)  }{2}.
\]
For $\alpha\in\mathbb{N}_{0}^{N}$ and $i\neq j$ by direct computation we
obtain:%
\begin{align}
\mathcal{B}_{ij}x^{\alpha} &  =\sum_{l=0}^{\alpha_{i}-\alpha_{j}}\left(
\frac{x_{j}}{x_{i}}\right)  ^{l}x^{\alpha},\text{ for }\alpha_{i}\geq
\alpha_{j},i<j,\label{bij1}\\
\mathcal{B}_{ij}x^{\alpha} &  =\sum_{l=0}^{\alpha_{i}-\alpha_{j}-1}\left(
\frac{x_{j}}{x_{i}}\right)  ^{l}x^{\alpha},\text{ for }\alpha_{i}\geq
\alpha_{j},i>j,\label{bij2}\\
\mathcal{B}_{ij}x^{\alpha} &  =-\sum_{l=1}^{\alpha_{j}-\alpha_{i}-1}\left(
\frac{x_{i}}{x_{j}}\right)  ^{l}x^{\alpha},\text{ for }\alpha_{i}<\alpha
_{j},i<j,\label{bij3}\\
\mathcal{B}_{ij}x^{\alpha} &  =-\sum_{l=1}^{\alpha_{j}-\alpha_{i}}\left(
\frac{x_{i}}{x_{j}}\right)  ^{l}x^{\alpha},\text{ for }\alpha_{i}<\alpha
_{j},i>j.\label{bij4}%
\end{align}
There is another invariant subspace structure for $\left\{  \mathcal{U}%
_{i}\left(  \kappa\right)  \right\}  $ besides the $\vartriangleright
$-triangular property. The purpose of the following arguments is to allow the
computation of certain coefficients of $\zeta_{\alpha}^{x}\left(
\kappa\right)  $ crucial in the arguments of Section 5.

\begin{definition}
For $1\leq s\leq N$ and $n\geq1$ let%
\begin{align*}
I_{s,n}^{\left(  N\right)  }  &  :=\left\{  \alpha\in\mathbb{N}_{0}^{N}%
:\alpha_{i}<n\text{ for }1\leq i\leq s\text{, }\alpha_{i}\leq n\text{ for
}s+1\leq i\leq N\right\}  ,\\
P_{s,n}^{\left(  N\right)  }  &  :=\mathrm{span}_{\mathbb{Q}\left(
\kappa\right)  }\left\{  x^{\alpha}:\alpha\in I_{s,n}^{\left(  N\right)
}\right\}  .
\end{align*}

\end{definition}

Note that each $I_{s,n}^{\left(  N\right)  }$ is finite and $P_{s,n}^{\left(
N\right)  }$ is the direct sum of its homogeneous subspaces $P_{s,n}^{\left(
N\right)  }\cap\mathcal{P}_{k},k\geq0$.

\begin{lemma}
\label{UE1}Suppose $1\leq s\leq N$ and $n\geq1$, then $\mathcal{U}_{i}\left(
\kappa\right)  P_{s,n}^{\left(  N\right)  }\subset P_{s,n}^{\left(  N\right)
}$ for $1\leq i\leq N$, and $\mathrm{span}_{\mathbb{Q}\left(  \kappa\right)
}\left\{  \zeta_{\alpha}^{x}\left(  \kappa\right)  :\alpha\in I_{s,n}^{\left(
N\right)  }\right\}  =P_{s,n}^{\left(  N\right)  }$.
\end{lemma}

\begin{proof}
Let $\alpha\in I_{s,n}^{\left(  N\right)  }$. It suffices to show
$\mathcal{B}_{ij}x^{\alpha}\in P_{s,n}^{\left(  N\right)  }$ for all $i,j$. By
formulae \ref{bij1}-\ref{bij4} this is obvious for $1\leq i,j\leq s$ or
$s+1\leq i,j\leq N$, or $\max\left(  \alpha_{i},\alpha_{j}\right)  <n$. Only
the two cases $1\leq i\leq s,\alpha_{j}=n$ (thus $j>s)$ and $1\leq j\leq
s,\alpha_{i}=n$ (with $i>s$) remain to be considered. Formulae \ref{bij3} and
\ref{bij2} respectively show that $\mathcal{B}_{ij}x^{\alpha}\in
P_{s,n}^{\left(  N\right)  }$. For any $\alpha\in\mathbb{N}_{0}^{N}$ the
eigenvector $\zeta_{\alpha}^{x}\left(  \kappa\right)  $ is contained in the
orbit of $x^{\alpha}$ under the algebra generated by $\left\{  \mathcal{U}%
_{i}\left(  \kappa\right)  \right\}  $ hence $x^{\alpha}\in P_{s,n}^{\left(
N\right)  }$ implies $\zeta_{\alpha}^{x}\left(  \kappa\right)  \in
P_{s,n}^{\left(  N\right)  }$. That the span of $\left\{  \zeta_{\alpha}%
^{x}\left(  \kappa\right)  \right\}  $ is all of $P_{s,n}^{\left(  N\right)
}$ follows easily (dimension argument, for example).
\end{proof}

\begin{definition}
\label{insrt}For a partition $\lambda$ and an integer $s$ with $1\leq s\leq N$
define the insertion operator $\iota\left(  s;\lambda\right)  :\mathbb{N}%
_{0}^{N}\rightarrow\mathbb{N}_{0}^{N+\ell\left(  \lambda\right)  }$ as
follows: for $\alpha\in\mathbb{N}_{0}^{N}$
\[
\left(  \iota\left(  s;\lambda\right)  \alpha\right)  _{i}=\left\{
\begin{array}
[c]{ll}%
\alpha_{i}, & 1\leq i\leq s\\
\lambda_{i-s}, & s<i\leq s+\ell\left(  \lambda\right) \\
\alpha_{i-\ell\left(  \lambda\right)  }, & s+\ell\left(  \lambda\right)
<i\leq N+\ell\left(  \lambda\right)  .
\end{array}
\right.
\]

\end{definition}

The definition is only interesting when $\alpha\in I_{s,n}^{\left(  N\right)
}$ where $n=\lambda_{\ell\left(  \lambda\right)  }$, in which case the
following rank equations hold: let $\beta=\iota\left(  s;\lambda\right)
\alpha$ and $k=\ell\left(  \lambda\right)  $, then $r\left(  \beta,i\right)
=r\left(  \alpha,i\right)  +k$ for $1\leq i\leq s$, $r\left(  \beta,i\right)
=r\left(  \alpha,i-k\right)  +k$ for $s+k<i\leq N+k$, and $r\left(
\beta,i\right)  =i-s$ for $s+1\leq i\leq s+k$.

\begin{theorem}
\label{coefs1}Suppose $\lambda$ is a partition, $1\leq s\leq N$, $\alpha
,\beta\in I_{s,n}^{\left(  N\right)  }$ where $n=\lambda_{\ell\left(
\lambda\right)  }$, and $\alpha\vartriangleright\beta$, then $\mathrm{coef}%
\left(  \zeta_{\iota\left(  s,\lambda\right)  \alpha}^{x}\left(
\kappa\right)  ,\iota\left(  s,\lambda\right)  \beta\right)  =\mathrm{coef}%
\left(  \zeta_{\alpha}^{x}\left(  \kappa\right)  ,\beta\right)  $.
\end{theorem}

It suffices to prove this for $\ell\left(  \lambda\right)  =1$ because then
one can insert one part of $\lambda$ at a time in nondecreasing order:
explicitly let $\lambda^{\left(  j\right)  }=\left(  \lambda_{\ell\left(
\lambda\right)  +1-j},\lambda_{\ell\left(  \lambda\right)  +2-j}%
,\ldots,\lambda_{\ell\left(  \lambda\right)  }\right)  $ for $1\leq j\leq
\ell\left(  \lambda\right)  $, then $\iota\left(  s,\left(  \lambda
_{\ell\left(  \lambda\right)  -j}\right)  \right)  \iota\left(  s,\lambda
^{\left(  j\right)  }\right)  =\iota\left(  s,\lambda^{\left(  j+1\right)
}\right)  $; also if $\alpha,\beta\in I_{s,n}^{\left(  N\right)  }$ then
$\iota\left(  s,\lambda^{\left(  j\right)  }\right)  \alpha\in I_{s,k}%
^{\left(  N+j\right)  }$ where $k=\lambda_{\ell\left(  \lambda\right)  +1-j}$
and $\alpha\vartriangleright\beta$ implies $\iota\left(  s,\lambda^{\left(
j\right)  }\right)  \alpha\vartriangleright\iota\left(  s,\lambda^{\left(
j\right)  }\right)  \beta$.

For arbitrary $M\geq1$ let $\mathcal{P}^{\left(  M\right)  }=\mathrm{span}%
_{\mathbb{Q}\left(  \kappa\right)  }\left\{  x^{\alpha}:\alpha\in
\mathbb{N}_{0}^{M}\right\}  $ and let $\mathcal{U}_{i}^{\left(  M\right)
}\left(  \kappa\right)  $ denote the operator $\mathcal{U}_{i}\left(
\kappa\right)  $ for $M$ variables. For $M>N$ let $\pi_{MN}$ be the projection
from $\mathcal{P}^{\left(  M\right)  }$ onto $\mathcal{P}^{\left(  N\right)
}$ defined by setting $x_{N+1}=x_{N+2}=\ldots=x_{M}=0$. The coefficients of
the $\zeta_{\alpha}^{x}$ do not depend on the number of variables (that is
$\mathrm{coef}\left(  \zeta_{\alpha}^{x},\beta\right)  $ is independent of
$N\geq\max\left(  \ell\left(  \alpha\right)  ,\ell\left(  \beta\right)
\right)  $) because of the intertwining relation
\begin{equation}
\pi_{MN}\mathcal{U}_{i}^{\left(  M\right)  }\left(  \kappa\right)  =\left(
\mathcal{U}_{i}^{\left(  N\right)  }\left(  \kappa\right)  +\left(
M-N\right)  \kappa\right)  \pi_{MN}, \label{pimn}%
\end{equation}
for $1\leq i\leq N<M$. Fix integers $n,s$ with $n\geq1$ and $1\leq s\leq N$
and define the map $\iota_{s,n}:\mathcal{P}^{\left(  N\right)  }%
\rightarrow\mathcal{P}^{\left(  N+1\right)  }$ by $\iota_{s,n}x^{\alpha
}=x^{\beta}$ for $\alpha\in\mathbb{N}_{0}^{N}$ and $\beta=\iota\left(
s,\left(  n\right)  \right)  \alpha=\left(  \alpha_{1},\ldots,\alpha
_{s},n,\alpha_{s+1},\ldots,\alpha_{N}\right)  $ and extending by linearity to
all polynomials. Direct computation yields the identities:%
\begin{align*}
\mathcal{U}_{i}^{\left(  N+1\right)  }\left(  \kappa\right)  \iota_{s,n}%
-\iota_{s,n}\mathcal{U}_{i}^{\left(  N\right)  }\left(  \kappa\right)   &
=\kappa\mathcal{B}_{i,s+1}\iota_{s,n},\text{ for }1\leq i\leq s,\\
\mathcal{U}_{i+1}^{\left(  N+1\right)  }\left(  \kappa\right)  \iota
_{s,n}-\iota_{s,n}\mathcal{U}_{i}^{\left(  N\right)  }\left(  \kappa\right)
&  =\kappa\mathcal{B}_{i+1,s+1}\iota_{s,n},\text{ for }s+1\leq i\leq N.
\end{align*}

We show that if $\alpha\in I_{s,n}^{\left(  N\right)  }$ then $\iota
_{s,n}\zeta_{\alpha}^{x}\left(  \kappa\right)  $ is congruent to $\zeta
_{\iota\left(  s,\left(  n\right)  \right)  \alpha}^{x}\left(  \kappa\right)
$ modulo the subspace $P_{s+1,n}^{\left(  N+1\right)  }$. To illustrate the
argument, suppose there is a linear operator $\mathcal{V}$ with an invariant
subspace $E$ and there is a vector $f$ and number $c$ so that $\mathcal{V}%
f-cf\in E$ then $f-\left(  \left(  \mathcal{V}-c\right)  |_{E}\right)
^{-1}\left(  \mathcal{V}f-cf\right)  $ is an eigenvector of $\mathcal{V}$ with
eigenvalue $c$, provided that the restriction of $\mathcal{V}-c$ to $E$ is
invertible. This can be adapted for simultaneous eigenvectors of pairwise
commuting operators by extending the base field $\mathbb{Q}\left(
\kappa\right)  $, adjoining another formal variable (transcendental) $v$ and
considering just one operator $\sum_{i=1}^{N}v^{i}\mathcal{U}_{i}^{\left(
N\right)  }\left(  \kappa\right)  $ (or $\sum_{i=1}^{N+1}v^{i}\mathcal{U}%
_{i}^{\left(  N+1\right)  }\left(  \kappa\right)  $, as appropriate). The
eigenvalues $\sum_{i=1}^{N}v^{i}\xi_{i}\left(  \alpha;\kappa\right)  $ are
simple ($\alpha\in\mathbb{N}_{0}^{N}$ and generic $\kappa$). Denote the field
$\mathbb{Q}\left(  \kappa,v\right)  $ by $\mathbb{K}.$

\begin{lemma}
Suppose $1\leq s\leq N$ and $n\geq1$. If $\alpha\in I_{s,n}^{\left(  N\right)
}$ then $\iota_{s,n}\zeta_{\alpha}^{x}\left(  \kappa\right)  =\zeta
_{\iota\left(  s,\left(  n\right)  \right)  \alpha}^{x}\left(  \kappa\right)
+f_{\alpha}$ for some $f_{\alpha}\in P_{s+1,n}^{\left(  N+1\right)  }$.
\end{lemma}

\begin{proof}
First we show $\mathcal{B}_{i,s+1}\iota_{s,n}P_{s,n}^{\left(  N\right)
}\subset P_{s+1,n}^{\left(  N+1\right)  }$ for $i\neq s+1$. Let $\alpha\in
I_{s,n}^{\left(  N\right)  }$ and $\beta=\iota\left(  s,\left(  n\right)
\right)  \alpha$. For $1\leq i\leq s$ by Formula \ref{bij3} $\mathcal{B}%
_{i,s+1}x^{\beta}=-\sum_{l=1}^{n-\alpha_{i}-1}\left(  \frac{x_{i}}{x_{s+1}%
}\right)  ^{l}x^{\beta}$ with the key (change from $\beta$) terms being
$x_{i}^{\alpha_{i}+l}x_{s+1}^{n-l}$ where $\alpha_{i}+1\leq\alpha_{i}+l\leq
n-1$ and $\alpha_{i}+1\leq n-l\leq n-1$; if $\alpha_{i}=n-1$ then
$\mathcal{B}_{i,s+1}x^{\beta}=0$. Suppose $s+2\leq i\leq N+1$; if
$\alpha_{i-1}=n=\beta_{i}$ then $\mathcal{B}_{i,s+1}x^{\beta}=0$ by Formula
\ref{bij2}, if $\alpha_{i-1}=\beta_{i}<n$ then by Formula \ref{bij4}
$\mathcal{B}_{i,s+1}x^{\beta}=-\sum_{l=1}^{n-\beta_{i}}\left(  \frac{x_{i}%
}{x_{s+1}}\right)  ^{l}x^{\beta}$ with key terms $x_{i}^{\beta_{i}+l}%
x_{s+1}^{n-l}$ where $\beta_{i}+1\leq\beta_{i}+l\leq n$ and $\beta_{i}\leq
n-l\leq n-1$. Thus $\mathcal{B}_{i,s+1}x^{\beta}\in P_{s+1,n}^{\left(
N+1\right)  }$.

Temporarily we use a superscript on the eigenvalues $\xi_{i}\left(
\alpha;\kappa\right)  $ to indicate the number of variables, then%
\begin{align*}
\xi_{i}^{\left(  N+1\right)  }\left(  \beta;\kappa\right)   &  =\left(
N+1-r\left(  \beta,i\right)  \right)  \kappa+\beta_{i}+1\\
&  =\left(  N+1-\left(  r\left(  \alpha,i\right)  +1\right)  \right)
\kappa+\alpha_{i}+1=\xi_{i}^{\left(  N\right)  }\left(  \alpha;\kappa\right)
\end{align*}
for $1\leq i\leq s$ and, similarly, $\xi_{i}^{\left(  N+1\right)  }\left(
\beta;\kappa\right)  =\xi_{i-1}^{\left(  N\right)  }\left(  \alpha
;\kappa\right)  $ for $s+2\leq i\leq N+1.$ The eigenvalues $\left\{  \xi
_{i}^{\left(  N+1\right)  }\left(  \beta;\kappa\right)  :1\leq i\leq N+1,i\neq
s+1\right\}  $ and the degree of homogeneity $\left\vert \beta\right\vert
=\left\vert \alpha\right\vert +n$ determine $\zeta_{\beta}^{x}\left(
\kappa\right)  $ uniquely, subject to $\mathrm{coef}\left(  \zeta_{\beta}%
^{x}\left(  \kappa\right)  ,\beta\right)  =1$, because $\sum_{i=1}%
^{N+1}\mathcal{U}_{i}^{\left(  N+1\right)  }\left(  \kappa\right)
=N+1+\sum_{i=1}^{N+1}x_{i}\frac{\partial}{\partial x_{i}}+\kappa\frac{N\left(
N+1\right)  }{2}$. Let%
\[
\mathcal{V}:=\left(  \sum_{i=1}^{s}+\sum_{i=s+2}^{N+1}\right)  v^{i}%
\mathcal{U}_{i}^{\left(  N+1\right)  }\left(  \kappa\right)  .
\]
The polynomials $\left\{  \zeta_{\gamma}^{x}\left(  \kappa\right)  :\gamma
\in\mathbb{N}_{0}^{N+1},\left\vert \gamma\right\vert =\left\vert
\alpha\right\vert +n\right\}  $ form a basis of eigenvectors of $\mathcal{V}$
for $E:=\mathrm{span}_{\mathbb{K}}\left\{  x^{\gamma}:\gamma\in\mathbb{N}%
_{0}^{N+1},\left\vert \gamma\right\vert =\left\vert \alpha\right\vert
+n\right\}  $ and each eigenvalue is simple. Let $F:=\mathrm{span}%
_{\mathbb{K}}\left\{  x^{\gamma}:\gamma\in I_{s+1,n}^{\left(  N+1\right)
},\left\vert \gamma\right\vert =\left\vert \alpha\right\vert +n\right\}  $
then $\mathcal{V}F\subset F$ by Lemma \ref{UE1} . Finally consider
\begin{align*}
\mathcal{V}\iota_{s,n}\zeta_{\alpha}^{x}\left(  \kappa\right)   &  =\sum
_{i=1}^{s}v^{i}\left(  \iota_{s,n}\mathcal{U}_{i}^{\left(  N\right)  }%
+\kappa\mathcal{B}_{i,s+1}\iota_{s,n}\right)  \zeta_{\alpha}^{x}\left(
\kappa\right)  +\\
&  +\sum_{i=s+2}^{N+1}v^{i}\left(  \iota_{s,n}\mathcal{U}_{i-1}^{\left(
N\right)  }+\kappa\mathcal{B}_{i,s+1}\iota_{s,n}\right)  \zeta_{\alpha}%
^{x}\left(  \kappa\right)  \\
&  =\sum_{i=1,i\neq s+1}^{N+1}v^{i}\xi_{i}^{\left(  N+1\right)  }\left(
\beta;i\right)  \iota_{s,n}\zeta_{\alpha}^{x}\left(  \kappa\right)
+h_{\alpha},
\end{align*}
where $h_{\alpha}=\kappa\left(  \sum_{i=1}^{s}+\sum_{i=s+2}^{N+1}\right)
v^{i}\mathcal{B}_{i,s+1}\iota_{s,n}\zeta_{\alpha}^{x}\left(  \kappa\right)  $
and $h_{\alpha}\in F$, since $\zeta_{\alpha}^{x}\left(  \kappa\right)  \in
P_{s,n}^{\left(  N\right)  }$. Let $\mathcal{V}_{\beta}$ be the restriction of
$\mathcal{V-}\sum_{i=1,i\neq s+1}^{N+1}v^{i}\xi_{i}^{\left(  N+1\right)
}\left(  \beta;i\right)  $ to the invariant subspace $F$ and let $f_{\alpha
}=\mathcal{V}_{\beta}^{-1}h_{\alpha}$, then $\iota_{s,n}\zeta_{\alpha}%
^{x}\left(  \kappa\right)  -f_{\alpha}=\zeta_{\beta}^{x}\left(  \kappa\right)
$ because $\mathrm{coef}\left(  f_{\alpha},\beta\right)  =0$ and
$\mathrm{coef}\left(  \iota_{s,n}\zeta_{\alpha}^{x}\left(  \kappa\right)
,\beta\right)  =\mathrm{coef}\left(  \zeta_{\alpha}^{x}\left(  \kappa\right)
,\alpha\right)  =1$. Since $f_{\alpha}=\iota_{s,n}\zeta_{\alpha}^{x}\left(
\kappa\right)  -\zeta_{\beta}^{x}\left(  \kappa\right)  $ the coefficients of
$f_{\alpha}$ are in $\mathbb{Q}\left(  \kappa\right)  $.
\end{proof}

\begin{corollary}
Suppose $\alpha,\gamma\in I_{s,n}^{\left(  N\right)  }$ then%
\[
\mathrm{coef}\left(  \zeta_{\iota\left(  s,\left(  n\right)  \right)  \alpha
}^{x}\left(  \kappa\right)  ,\iota\left(  s,\left(  n\right)  \right)
\gamma\right)  =\mathrm{coef}\left(  \zeta_{\alpha}^{x}\left(  \kappa\right)
,\gamma\right)  .
\]

\end{corollary}

\begin{proof}
By definition $\mathrm{coef}\left(  \zeta_{\alpha}^{x}\left(  \kappa\right)
,\gamma\right)  =\mathrm{coef}\left(  \iota_{s,n}\zeta_{\alpha}^{x}\left(
\kappa\right)  ,\iota\left(  s,\left(  n\right)  \right)  \gamma\right)  $.
Also \newline$\left(  \iota\left(  s,\left(  n\right)  \right)  \gamma\right)
_{s+1}=n$ and thus $\mathrm{coef}\left(  f,\iota\left(  s,\left(  n\right)
\right)  \gamma\right)  =0$ for any $f\in P_{s+1,n}^{\left(  N+1\right)  }$.
\end{proof}

This completes the proof of Theorem \ref{coefs1}.

The poles of the coefficients of $\zeta_{\alpha}^{x}\left(  \kappa\right)  $
play a key role in the analysis of singular polynomials. Knop and Sahi
\cite{KS} found an algorithm for the evaluation of the coefficients. It uses
the idea of extending the definition of Ferrers diagrams to compositions and
associating a hook-length to each node in the diagram. The \textit{Ferrers
diagram} of a composition $\alpha\in\mathbb{N}_{0}^{N}$ is the set $\left\{
\left(  i,j\right)  :1\leq i\leq\ell\left(  \alpha\right)  ,0\leq j\leq
\alpha_{i}\right\}  .$ For each node $\left(  i,j\right)  $ with $1\leq
j\leq\alpha_{i}$ \ there are two special subsets of the Ferrers diagram, the
\textit{arm} $\left\{  \left(  i,l\right)  :j<l\leq\alpha_{i}\right\}  $ and
the \textit{leg} $\left\{  \left(  l,j\right)  :l>i,j\leq\alpha_{l}\leq
\alpha_{i}\right\}  \cup\left\{  \left(  l,j-1\right)  :l<i,j-1\leq\alpha
_{l}<\alpha_{i}\right\}  $. The node itself, the arm and the leg make up the
\textit{hook}. The definition of hooks for compositions is from \cite[p.15]%
{KS}. The cardinality of the leg is called the leg-length, formalized by the following:

\begin{definition}
For $\alpha\in\mathbb{N}_{0}^{N},1\leq i\leq\ell\left(  \alpha\right)  $ and
$1\leq j\leq\alpha_{i}$ the leg-length is%
\begin{align*}
L\left(  \alpha;i,j\right)   &  :=\#\left\{  l:l>i,j\leq\alpha_{l}\leq
\alpha_{i}\right\} \\
&  +\#\left\{  l:l<i,j\leq\alpha_{l}+1\leq\alpha_{i}\right\}  .
\end{align*}
For $t\in\mathbb{Q}\left(  \kappa\right)  $ the \textit{hook-length }and the
hook-length product for $\alpha$ are given by
\begin{align*}
h\left(  \alpha,t;i,j\right)   &  =\left(  \alpha_{i}-j+t+\kappa L\left(
\alpha;i,j\right)  \right) \\
h\left(  \alpha,t\right)   &  =\prod_{i=1}^{\ell\left(  \alpha\right)  }%
\prod_{j=1}^{\alpha_{i}}h\left(  \alpha,t;i,j\right)  ,
\end{align*}

\end{definition}

Note that the indices $\left\{  i:\alpha_{i}=0\right\}  $ are omitted in the
product $h\left(  \alpha,t\right)  $. In \cite{D2} and \cite{DX} we used the
notation
\[
\mathcal{E}_{\varepsilon}\left(  \alpha\right)  =\prod\left\{  1+\frac
{\varepsilon\kappa}{\kappa\left(  r\left(  \alpha,i\right)  -r\left(
\alpha,j\right)  \right)  +\alpha_{j}-\alpha_{i}}:i<j,\,\alpha_{i}<\alpha
_{j}\right\}  ,\varepsilon=\pm.
\]
The denominator also equals $\xi_{j}\left(  \alpha;\kappa\right)  -\xi
_{i}\left(  \alpha;\kappa\right)  $. The relation to $h\left(  \alpha
,t\right)  $ (for the values $t=1,\kappa+1$ which are of concern here) is the following:

\begin{lemma}
For $\alpha\in\mathbb{N}_{0}^{N}$, $h\left(  \alpha,\kappa+1\right)  =h\left(
\alpha^{+},\kappa+1\right)  \mathcal{E}_{+}\left(  \alpha\right)  $ and
$h\left(  \alpha,1\right)  =\dfrac{h\left(  \alpha^{+},1\right)  }%
{\mathcal{E}_{-}\left(  \alpha\right)  }$.
\end{lemma}

\begin{proof}
We use induction on adjacent transpositions. The statements are true for
$\alpha=\alpha^{+}$. Fix $\alpha^{+}$ and suppose $\alpha_{i}>\alpha_{i+1}$
for some $i$. Let $\sigma=\left(  i,i+1\right)  $. Consider the ratio
$\dfrac{h\left(  \sigma\alpha,t\right)  }{h\left(  \alpha,t\right)  }$. The
only node whose hook-length changes (in the sense of interchanging rows $i$
and $i+1$ of the Ferrers diagram) is $\left(  i,\alpha_{i+1}+1\right)  $.
Explicitly $h\left(  \sigma\alpha,t;s,j\right)  =h\left(  \alpha,t;s,j\right)
$ for $s\neq i,i+1$ and $1\leq j\leq\alpha_{s}$, $h\left(  \sigma
\alpha,t;i,j\right)  =h\left(  \alpha,t;i+1,j\right)  $ for $1\leq j\leq
\alpha_{i+1}$ and $h\left(  \sigma\alpha,t;i+1,j\right)  =h\left(
\alpha,t;i,j\right)  $ for $1\leq j\leq\alpha_{i}$ except for $j=\alpha
_{i+1}+1$. Thus $\dfrac{h\left(  \sigma\alpha,t\right)  }{h\left(
\alpha,t\right)  }=\dfrac{h\left(  \sigma\alpha,t;i+1,\alpha_{i+1}+1\right)
}{h\left(  \alpha,t;i,\alpha_{i+1}+1\right)  }$. Note that $L\left(
\sigma\alpha;i+1,\alpha_{i+1}+1\right)  =L\left(  \alpha;i,\alpha
_{i+1}+1\right)  +1$ (the node $\left(  i,\alpha_{i+1}\right)  $ is adjoined
to the leg). Let%
\begin{align*}
E_{1} &  =\left\{  s:s\leq i,\alpha_{s}\geq\alpha_{i}\right\}  \cup\left\{
s:s>i,\alpha_{s}>\alpha_{i}\right\}  ,\\
E_{2} &  =\left\{  s:s\leq i+1,\alpha_{s}\geq\alpha_{i+1}\right\}
\cup\left\{  s:s>i+1,\alpha_{s}>\alpha_{i+1}\right\}  ,
\end{align*}
thus by definition $r\left(  \alpha,i\right)  =\#E_{1}$ and $r\left(
\alpha,i+1\right)  =\#E_{2}$. Now $E_{1}\subset E_{2}$ thus $r\left(
\alpha,i+1\right)  -r\left(  \alpha,i\right)  =\#\left(  E_{2}\backslash
E_{1}\right)  $ and $E_{2}\backslash E_{1}=\left\{  s:s<i,\alpha_{i}%
>\alpha_{s}\geq\alpha_{i+1}\right\}  \cup\left\{  i\right\}  \cup\left\{
s:s>i+1,\alpha_{i}\geq\alpha_{s}>\alpha_{i+1}\right\}  $. This shows that
$\#\left(  E_{2}\backslash E_{1}\right)  =1+L\left(  \alpha;i,\alpha
_{i+1}+1\right)  ,$ and
\begin{align*}
h\left(  \alpha,t;i,\alpha_{i+1}+1\right)   &  =\kappa\left(  r\left(
\alpha,i+1\right)  -r\left(  \alpha,i\right)  -1\right)  +t+\alpha_{i}%
-\alpha_{i+1}-1,\\
h\left(  \sigma\alpha,t;i+1,\alpha_{i+1}+1\right)   &  =\kappa\left(  r\left(
\alpha,i+1\right)  -r\left(  \alpha,i\right)  \right)  +t+\alpha_{i}%
-\alpha_{i+1}-1.
\end{align*}
Thus
\begin{align*}
\dfrac{h\left(  \sigma\alpha,\kappa+1;i+1,\alpha_{i+1}+1\right)  }{h\left(
\alpha,\kappa+1;i,\alpha_{i+1}+1\right)  } &  =\dfrac{\kappa\left(  r\left(
\alpha,i+1\right)  -r\left(  \alpha,i\right)  +1\right)  +\alpha_{i}%
-\alpha_{i+1}}{\kappa\left(  r\left(  \alpha,i+1\right)  -r\left(
\alpha,i\right)  \right)  +\alpha_{i}-\alpha_{i+1}}\\
&  =1+\frac{\kappa}{\kappa\left(  r\left(  \alpha,i+1\right)  -r\left(
\alpha,i\right)  \right)  +\alpha_{i}-\alpha_{i+1}}\\
&  =\mathcal{E}_{+}\left(  \sigma\alpha\right)  /\mathcal{E}_{+}\left(
\alpha\right)  ;
\end{align*}
the latter equation is proven in Theorem 8.5.8,from \cite[p.302]{DX}, and
\begin{align*}
\dfrac{h\left(  \sigma\alpha,1;i+1,\alpha_{i+1}+1\right)  }{h\left(
\alpha,1;i,\alpha_{i+1}+1\right)  } &  =\dfrac{\kappa\left(  r\left(
\alpha,i+1\right)  -r\left(  \alpha,i\right)  \right)  +\alpha_{i}%
-\alpha_{i+1}}{\kappa\left(  r\left(  \alpha,i+1\right)  -r\left(
\alpha,i\right)  -1\right)  +\alpha_{i}-\alpha_{i+1}}\\
&  =\left(  1-\frac{\kappa}{\kappa\left(  r\left(  \alpha,i+1\right)
-r\left(  \alpha,i\right)  \right)  +\alpha_{i}-\alpha_{i+1}}\right)  ^{-1}\\
&  =\mathcal{E}_{-}\left(  \alpha\right)  /\mathcal{E}_{-}\left(  \sigma
\alpha\right)  .
\end{align*}
Thus $h\left(  \alpha,\kappa+1\right)  $ and $h\left(  \alpha^{+}%
,\kappa+1\right)  \mathcal{E}_{+}\left(  \alpha\right)  $ have the same
transformation properties under adjacent transpositions and hence are equal.
Similarly $h\left(  \alpha,1\right)  =\frac{h\left(  \alpha^{+},1\right)
}{\mathcal{E}_{-}\left(  \alpha\right)  }$.
\end{proof}

Knop and Sahi \cite[Theorem 5.1]{KS} showed that $h\left(  \alpha
,\kappa+1\right)  \zeta_{\alpha}^{x}\left(  \kappa\right)  $ has all
coefficients in $\mathbb{N}_{0}\left[  \kappa\right]  $ for each $\alpha
\in\mathbb{N}_{0}^{N}$. When $\kappa$ takes on a negative rational number
$\kappa_{0}$ it may happen that two different compositions have the same
eigenvalues $\left(  \xi_{i}\left(  \alpha;\kappa_{0}\right)  \right)
_{i=1}^{N}$ so one can not claim the existence of a basis of simultaneous
eigenvectors of $\left\{  \mathcal{U}_{i}\left(  \kappa_{0}\right)  :1\leq
i\leq N\right\}  $. We recall the following from \cite{D2}.

\begin{definition}
Let $\alpha,\beta\in\mathbb{N}_{0}^{N}$ and let $m,n\in\mathbb{N}$ with
$\gcd\left(  m,n\right)  =1$ then say $\left(  \alpha,\beta\right)  $ is a
$\left(  -\frac{m}{n}\right)  $-critical pair (for $\alpha$) if $\alpha
\vartriangleright\beta$ and $\left(  n\kappa+m\right)  $ divides $\left(
r\left(  \beta,i\right)  -r\left(  \alpha,i\right)  \right)  \kappa+\alpha
_{i}-\beta_{i}$ (in $\mathbb{Q}\left[  \kappa\right]  $) for $1\leq i\leq N$.
\end{definition}

The definition implies $\xi_{i}\left(  \alpha;-\frac{m}{n}\right)  =\xi
_{i}\left(  \beta;-\frac{m}{n}\right)  $ for $1\leq i\leq N$. We can deduce
the existence of simple poles at $\kappa=-\frac{m}{n}$ in a certain coefficient.

\begin{lemma}
\label{kpole}Suppose $\alpha,\beta\in\mathbb{N}_{0}^{N}$, $h\left(
\alpha,\kappa+1\right)  $ has a simple zero at $\kappa_{0}\in\mathbb{Q}$ and
$\left(  \alpha,\beta\right)  $ is the unique $\kappa_{0}$-critical pair for
$\alpha$, then $\mathrm{coef}\left(  \zeta_{\alpha}^{x}\left(  \kappa\right)
,\beta\right)  $ has a simple pole at $\kappa_{0}$.
\end{lemma}

\begin{proof}
Since $\mathrm{coef}\left(  \zeta_{\alpha}^{x}\left(  \kappa\right)
,\beta\right)  $ is independent of the number of variables $N$ provided
$N\geq\max\left(  \ell\left(  \alpha\right)  ,\ell\left(  \beta\right)
\right)  $ we may assume $N=\ell\left(  \alpha\right)  +\left\vert
\alpha\right\vert $. Let $\gamma=\left(  0^{\ell\left(  \alpha\right)
},1^{\left\vert \alpha\right\vert }\right)  \in\mathbb{N}_{0}^{N}$ then by
\cite{KS} $\mathrm{coef}\left(  \zeta_{\alpha}^{x}\left(  \kappa\right)
,\gamma\right)  =\left(  \left\vert \alpha\right\vert \right)  !\kappa
^{\left\vert \alpha\right\vert }/h\left(  \alpha,\kappa+1\right)  $. Let
$f=\lim\limits_{\kappa\rightarrow\kappa_{0}}\left(  \kappa-\kappa_{0}\right)
\zeta_{\alpha}^{x}\left(  \kappa\right)  $ which exists as a polynomial over
$\mathbb{Q}$ by hypothesis and is not zero because\newline$\lim\limits_{\kappa
\rightarrow\kappa_{0}}\left(  \kappa-\kappa_{0}\right)  \mathrm{coef}\left(
\zeta_{\alpha}^{x}\left(  \kappa\right)  ,\gamma\right)  \neq0$. The
polynomial $f$ is a simultaneous eigenvector for $\left\{  \mathcal{U}%
_{i}\left(  \kappa_{0}\right)  :1\leq i\leq N\right\}  $ because
$\mathcal{U}_{i}\left(  \kappa_{0}\right)  f=\lim\limits_{\kappa
\rightarrow\kappa_{0}}\left(  \kappa-\kappa_{0}\right)  \mathcal{U}_{i}\left(
\kappa\right)  \zeta_{\alpha}^{x}\left(  \kappa\right)  =\xi_{i}\left(
\alpha;\kappa_{0}\right)  f$. Let $\gamma$ be a $\vartriangleright$-maximal
element of $\left\{  \delta\in\mathbb{N}_{0}^{N}:\mathrm{coef}\left(
f,\delta\right)  \neq0\right\}  $. By $\vartriangleright$-triangularity
$\xi_{i}\left(  \alpha;\kappa_{0}\right)  =\xi_{i}\left(  \delta;\kappa
_{0}\right)  $ for each $i$, thus $\delta=\alpha$ or $\delta=\beta$ by
definition of critical pairs. It is impossible for $\delta=\alpha$ since
$\mathrm{coef}\left(  f,\alpha\right)  =\lim\limits_{\kappa\rightarrow
\kappa_{0}}\left(  \kappa-\kappa_{0}\right)  =0$ hence $\delta=\beta$. So
$\mathrm{coef}\left(  f,\beta\right)  =\lim\limits_{\kappa\rightarrow
\kappa_{0}}\left(  \kappa-\kappa_{0}\right)  \mathrm{coef}\left(
\zeta_{\alpha}^{x}\left(  \kappa\right)  ,\beta\right)  \neq0$.
\end{proof}

In the next sections the Lemma will be combined with Theorem \ref{coefs1}.

\begin{example}
The conceptual proof of the Lemma may be the only reasonably effective method.
For example in the next section we need the conclusion of Lemma \ref{kpole}
for $\mathrm{coef}\left(  \zeta_{\left(  5,6\right)  }^{x}\left(
\kappa\right)  ,\left(  2,0,3,3,3\right)  \right)  $, which arises for
$N=5,\tau=\left(  3,2\right)  ,\kappa_{0}=-\frac{3}{2}$. There is a
combinatorial formula for the coefficients of $\zeta_{\alpha}^{x}\left(
\kappa\right)  $ due to Knop and Sahi \cite{KS}, which requires a sum over
$3!\times1721$ configurations for this example (the factorial comes from
permuting the indices $\left(  3,4,5\right)  $). By direct (computer algebra)
calculations this coefficient equals%
\[
\frac{30\kappa^{3}(1+\kappa)^{2}(62\kappa^{3}+135\kappa^{2}+78\kappa
+40)}{(2\kappa+3)(2\kappa+5)(\kappa+2)^{2}(\kappa+3)^{2}(\kappa+4)(\kappa
+5)}.
\]
The expression suggests that there is no practical closed form.
\end{example}

We address the problem of the relationship of a simultaneous eigenvector of
$\left\{  \mathcal{U}_{i}\left(  \kappa_{0}\right)  :1\leq i\leq N\right\}  $
to the nonsymmetric Jack polynomials; namely, how can such a polynomial be
expressed as a limit as $\kappa\rightarrow\kappa_{0}$? For a given $\alpha
\in\mathbb{N}_{0}^{N}$ and $\kappa_{0}=-\frac{m}{n}$ let $C\left(
\alpha,\kappa_{0}\right)  =\left\{  \beta:\left(  \alpha,\beta\right)  \text{
is a }\kappa_{0}\text{-critical pair}\right\}  $. In the proof we again use
the field $\mathbb{K}=\mathbb{Q}\left(  \kappa,v\right)  $ and the operator
$\sum_{i=1}^{N}v^{i}\mathcal{U}_{i}\left(  \kappa\right)  $; otherwise to each
$\gamma\in E$ one has to associate some $i$ for which $\xi_{i}\left(
\alpha;\kappa_{0}\right)  \neq\xi_{i}\left(  \gamma;\kappa_{0}\right)  $. The
expressions we consider are all rational in $\kappa$ (now with values in
$\mathbb{Q}\left(  v\right)  $) so having no pole at $\kappa_{0}$ is
equivalent to being analytic in a neighborhood of $\kappa_{0}$.

\begin{theorem}
\label{eigvect0}Suppose for some $\alpha\in\mathbb{N}_{0}^{N}$ and $\kappa
_{0}=-\frac{m}{n}$ that there exists a simultaneous eigenvector $g_{\alpha
}=x^{\alpha}+\sum_{\beta\vartriangleleft\alpha}A_{\beta}x^{\beta}$ of
$\left\{  \mathcal{U}_{i}\left(  \kappa_{0}\right)  :1\leq i\leq N\right\}  $,
with the coefficients $A_{\beta}\in\mathbb{Q}$, then there are coefficients
$B_{\beta}\left(  \kappa\right)  \in\mathbb{K}$ defined for $\beta\in C\left(
\alpha,\kappa_{0}\right)  $ such that the polynomial $q_{\alpha}\left(
\kappa\right)  =\zeta_{\alpha}^{x}\left(  \kappa\right)  +\sum\left\{
B_{\beta}\left(  \kappa\right)  \zeta_{\beta}^{x}\left(  \kappa\right)
:\beta\in C\left(  \alpha,\kappa_{0}\right)  \right\}  $ has no pole at
$\kappa_{0}$ and $\lim\limits_{\kappa\rightarrow\kappa_{0}}q_{\alpha}\left(
\kappa\right)  =g_{\alpha}$.
\end{theorem}

\begin{proof}
By the triangularity property $\mathcal{U}_{i}\left(  \kappa_{0}\right)
g_{\alpha}=\xi_{i}\left(  \alpha;\kappa_{0}\right)  g_{\alpha}$. For generic
$\kappa$ there are coefficients $B_{\gamma}^{\prime}\left(  \kappa\right)  $
defined for all $\gamma\vartriangleleft\alpha$ so that $g_{\alpha}%
=\zeta_{\alpha}^{x}\left(  \kappa\right)  +\sum_{\gamma\vartriangleleft\alpha
}B_{\gamma}^{\prime}\left(  \kappa\right)  \zeta_{\gamma}^{x}\left(
\kappa\right)  $ (because the nonsymmetric Jack polynomials form a basis and
the change-of-basis matrix is unimodular and triangular). Let $E=\left\{
\gamma:\gamma\vartriangleleft\alpha,\gamma\notin C\left(  \alpha,\kappa
_{0}\right)  \right\}  $. Apply the operator%
\[
\mathcal{V}\left(  \kappa\right)  :=\prod\nolimits_{\gamma\in E}\frac
{\sum_{i=1}^{N}v^{i}\left(  \mathcal{U}_{i}\left(  \kappa\right)  -\xi
_{i}\left(  \gamma;\kappa\right)  \right)  }{\sum_{i=1}^{N}v^{i}\left(
\xi_{i}\left(  \alpha;\kappa\right)  -\xi_{i}\left(  \gamma;\kappa\right)
\right)  }%
\]
to both sides of the equation for $g_{\alpha}$, thus annihilating all
$\zeta_{\gamma}^{x}\left(  \kappa\right)  $ with $\gamma\in E$. The right hand
side becomes
\begin{align*}
\mathcal{V}\left(  \kappa\right)  g_{\alpha}  &  =\zeta_{\alpha}^{x}\left(
\kappa\right)  +\sum_{\beta\in C\left(  \alpha,\kappa_{0}\right)  }B_{\beta
}^{\prime}\left(  \kappa\right)  \left(  \prod\nolimits_{\gamma\in E}%
\frac{\sum_{i=1}^{N}v^{i}\left(  \xi_{i}\left(  \beta;\kappa\right)  -\xi
_{i}\left(  \gamma;\kappa\right)  \right)  }{\sum_{i=1}^{N}v^{i}\left(
\xi_{i}\left(  \alpha;\kappa\right)  -\xi_{i}\left(  \gamma;\kappa\right)
\right)  }\right)  \zeta_{\beta}^{x}\left(  \kappa\right) \\
&  =\zeta_{\alpha}^{x}\left(  \kappa\right)  +\sum_{\beta\in C\left(
\alpha,\kappa_{0}\right)  }B_{\beta}\left(  \kappa\right)  \zeta_{\beta}%
^{x}\left(  \kappa\right)  ,
\end{align*}
with the last equation implicitly defining the coefficients $B_{\beta}\left(
\kappa\right)  $.We use the operators $\mathcal{B}_{i}$ from Definition
\ref{defbij}. To evaluate $\mathcal{V}\left(  \kappa\right)  g_{\alpha}$
directly we consider%
\begin{align*}
&  \left(  \mathcal{U}_{i}\left(  \kappa\right)  -\xi_{i}\left(  \gamma
;\kappa\right)  \right)  g_{\alpha}-\left(  \xi_{i}\left(  \alpha
;\kappa\right)  -\xi_{i}\left(  \gamma;\kappa\right)  \right)  g_{\alpha}\\
&  =\left(  \mathcal{U}_{i}\left(  \kappa_{0}\right)  +\left(  \kappa
-\kappa_{0}\right)  \mathcal{B}_{i}-\xi_{i}\left(  \alpha;\kappa\right)
\right)  g_{\alpha}\\
&  =\left(  \xi_{i}\left(  \alpha;\kappa_{0}\right)  -\xi_{i}\left(
\alpha;\kappa\right)  +\left(  \kappa-\kappa_{0}\right)  \mathcal{B}%
_{i}\right)  g_{\alpha}\\
&  =\left(  \kappa-\kappa_{0}\right)  \left(  r\left(  \alpha,i\right)
-N+\mathcal{B}_{i}\right)  g_{\alpha}.
\end{align*}
Thus for each $\gamma\in E$ we have%
\[
\dfrac{\sum_{i=1}^{N}v^{i}\left(  \mathcal{U}_{i}\left(  \kappa\right)
-\xi_{i}\left(  \gamma;\kappa\right)  \right)  }{\sum_{i=1}^{N}v^{i}\left(
\xi_{i}\left(  \alpha;\kappa\right)  -\xi_{i}\left(  \gamma;\kappa\right)
\right)  }g_{\alpha}=g_{\alpha}+\left(  \kappa-\kappa_{0}\right)  \dfrac
{\sum_{i=1}^{N}v^{i}\left(  r\left(  \alpha,i\right)  -N+\mathcal{B}%
_{i}\right)  }{\sum_{i=1}^{N}v^{i}\left(  \xi_{i}\left(  \alpha;\kappa\right)
-\xi_{i}\left(  \gamma;\kappa\right)  \right)  }g_{\alpha}.
\]
The latter term has no pole at $\kappa_{0}$, since $\left(  \alpha
,\gamma\right)  $ is not a $\kappa_{0}$-critical pair. Apply this computation
repeatedly to obtain $\mathcal{V}\left(  \kappa\right)  g_{\alpha}=g_{\alpha
}+\left(  \kappa-\kappa_{0}\right)  p\left(  \kappa\right)  $, where $p$ is
polynomial in $x$, rational in $\kappa$ and has no pole at $\kappa_{0}$. Hence
set $q_{\alpha}\left(  \kappa\right)  =\mathcal{V}\left(  \kappa\right)
g_{\alpha}$, then $\lim\limits_{\kappa\rightarrow\kappa_{0}}q_{\alpha}\left(
\kappa\right)  =g_{\alpha}$, and this completes the proof.
\end{proof}

If we apply this result to the hypothetical singular polynomial described in
Theorem \ref{glambda}, that is $g_{\lambda}=x^{\lambda}+\sum_{\beta
\vartriangleleft\lambda}A_{\beta}x^{\beta},$ we obtain%
\[
\mathcal{V}\left(  \kappa\right)  g_{\lambda}=\zeta_{\lambda}^{x}\left(
\kappa\right)  +\sum\left\{  B_{\beta}\left(  \kappa\right)  \zeta_{\beta}%
^{x}\left(  \kappa\right)  :\beta\in C\left(  \lambda,\kappa_{0}\right)
\right\}  ,
\]
which has no pole at $\kappa_{0}$. More importantly, since $\mathcal{D}%
_{i}\left(  \kappa\right)  $ is polynomial in $\kappa$, the relation
$\lim\limits_{\kappa\rightarrow\kappa_{0}}\mathcal{D}_{i}\left(
\kappa\right)  \mathcal{V}\left(  \kappa\right)  g_{\lambda}=\mathcal{D}%
_{i}\left(  \kappa_{0}\right)  g_{\lambda}=0$ holds. This is a key ingredient
in the proof that $\tau_{2}<n$, because we can now apply the known formulae
for $\mathcal{D}_{i}\left(  \kappa\right)  \zeta_{\beta}^{x}\left(
\kappa\right)  $.

The basic step is the formula for $\mathcal{D}_{\ell\left(  \alpha\right)
}\left(  \kappa\right)  \zeta_{\alpha}^{x}\left(  \kappa\right)  $ for
$\alpha\in\mathbb{N}_{0}^{N}.$ The computation involves a cyclic shift. For
$1\leq i\leq N$ let $\varepsilon\left(  i\right)  \in\mathbb{N}_{0}^{N}$
denote the standard basis element, that is, $\varepsilon\left(  i\right)
_{j}=\delta_{ij}$.

\begin{definition}
\label{cyclic}For $1<k\leq N$ let $\theta_{k}=\left(  1,2\right)  \left(
2,3\right)  \ldots\left(  k-1,k\right)  \in S_{N}$, (thus, $\theta_{k}%
\alpha=\left(  \alpha_{k},\alpha_{1},\ldots,\alpha_{k-1},\alpha_{k+1}%
,\ldots\right)  $ for $\alpha\in\mathbb{N}_{0}^{N}$). If $\alpha\in
\mathbb{N}_{0}^{N}$ satisfies $\ell\left(  \alpha\right)  =k$ for $1<k\leq N$
set $\widetilde{\alpha}=\theta_{k}\left(  \alpha-\varepsilon\left(  k\right)
\right)  =\left(  \alpha_{k}-1,\alpha_{1},\ldots,\alpha_{k-1},0,\ldots\right)
$.
\end{definition}

In \cite{D2} the formula for $\mathcal{D}_{k}\left(  \kappa\right)  $ is
stated for the $p$-basis $\left\{  \zeta_{\alpha}\left(  \kappa\right)
\right\}  $. To use the result here it suffices to invoke the transformation
formula $\zeta_{\alpha}\left(  \kappa\right)  =\frac{h\left(  \alpha
,\kappa+1\right)  }{h\left(  \alpha,1\right)  }\zeta_{\alpha}^{x}\left(
\kappa\right)  $ for $\alpha\in\mathbb{N}_{0}^{N}$. The ratio does not have to
be computed explicitly since only the values of $\frac{h\left(  \alpha
,t\right)  }{h\left(  \widetilde{\alpha},t\right)  }$ for $t=1,\kappa+1$ are needed.

\begin{lemma}
\label{hookd}Let $\alpha\in\mathbb{N}_{0}^{N}$ and suppose $\ell\left(
\alpha\right)  =k$ then $\dfrac{h\left(  \alpha,t\right)  }{h\left(
\widetilde{\alpha},t\right)  }=\left(  k-r\left(  \alpha,k\right)  \right)
\kappa+t+\alpha_{k}-1$.
\end{lemma}

\begin{proof}
Heuristically the Ferrers diagram for $\widetilde{\alpha}$ is produced from
that of $\alpha$ by deleting the node at $\left(  k,1\right)  $ and moving the
remainder of row $k$ to the top (row zero); then every node still has the same
hook-length and the required ratio is $h\left(  a,t;k,1\right)  $. Explicitly,
$h\left(  \widetilde{\alpha},t;i,j\right)  =h\left(  \alpha,t;i-1,j\right)  $
for $2\leq i\leq k,1\leq j\leq\alpha_{i-1}$ and $h\left(  \widetilde{\alpha
},t;1,j-1\right)  =h\left(  \alpha,t;k,j\right)  $ for $2\leq j\leq\alpha_{k}$
because%
\begin{align*}
L\left(  \alpha;k,j\right)   &  =\#l:l<k,j\leq\alpha_{l}+1\leq\alpha_{k}\\
&  =\#\left\{  l:1<l,j-1\leq\alpha_{l-1}\leq\alpha_{k}-1\right\}  =L\left(
\widetilde{\alpha};1,j-1\right)  .
\end{align*}
Also%
\begin{align*}
L\left(  \alpha;k,1\right)   &  =\#\left\{  l:l<k,1\leq\alpha_{l}+1\leq
\alpha_{k}\right\}  =\#\left\{  l:l<k,0\leq\alpha_{l}<\alpha_{k}\right\} \\
&  =k-\#\left\{  l:l\leq k,\alpha_{l}\geq\alpha_{k}\right\}  =k-r\left(
\alpha,k\right)  ,
\end{align*}
thus $h\left(  \alpha,t;1,k\right)  =\kappa\left(  k-r\left(  \alpha,k\right)
\right)  +t+\alpha_{k}-1$.
\end{proof}

\begin{proposition}
\label{zdiff}Let $\alpha\in\mathbb{N}_{0}^{N}$ and suppose $\ell\left(
\alpha\right)  =k$ then%
\[
\mathcal{D}_{k}\left(  \kappa\right)  \zeta_{\alpha}^{x}\left(  \kappa\right)
=\dfrac{\left(  k-r\left(  \alpha,k\right)  \right)  \kappa+\alpha_{k}%
}{\left(  k+1-r\left(  \alpha,k\right)  \right)  \kappa+\alpha_{k}}\left(
\left(  N+1-r\left(  \alpha,k\right)  \right)  \kappa+\alpha_{k}\right)
\theta_{k}^{-1}\zeta_{\widetilde{\alpha}}^{x}\left(  \kappa\right)  .
\]

\end{proposition}

\begin{proof}
In \cite[Theorem 3.5]{D2} it was shown that%
\[
\mathcal{D}_{k}\left(  \kappa\right)  \zeta_{\alpha}\left(  \kappa\right)
=\left(  \left(  N+1-r\left(  \alpha,k\right)  \right)  \kappa+\alpha
_{k}\right)  \theta_{k}^{-1}\zeta_{\widetilde{\alpha}}\left(  \kappa\right)
.
\]
To modify this equation to hold for the $x$-monic polynomials multiply the
right hand side by $\dfrac{h\left(  \alpha,1\right)  h\left(  \widetilde
{\alpha},\kappa+1\right)  }{h\left(  \widetilde{\alpha},1\right)  h\left(
\alpha,\kappa+1\right)  }=\dfrac{\left(  k-r\left(  \alpha,k\right)  \right)
\kappa+\alpha_{k}}{\left(  k+1-r\left(  \alpha,k\right)  \right)
\kappa+\alpha_{k}}$.
\end{proof}

The last topic for the section is the action of $\mathcal{D}_{i}\left(
\kappa\right)  $ with respect to the order $\vartriangleright$. In the lemma
the operator is modified to be degree-preserving to simplify the statement.

\begin{lemma}
\label{dcoefs}Suppose $\alpha\in\mathbb{N}_{0}^{N}$ and $1\leq i\leq N$, if
$\mathrm{coef}\left(  x_{i}\mathcal{D}_{i}\left(  \kappa\right)  x^{\beta
},\alpha\right)  \neq0$ then $\beta=\alpha$ or $\beta^{+}\succ\alpha^{+}$ or
$\beta=\left(  i,j\right)  \alpha$ with $\alpha_{i}>\alpha_{j},1\leq j\leq N$.
\end{lemma}

\begin{proof}
By direct computation for $\beta\in\mathbb{N}_{0}^{N}$ we have
\[
x_{i}\mathcal{D}_{i}\left(  \kappa\right)  x^{\beta}=\beta_{i}x^{\beta}%
+\kappa\sum_{\beta_{j}<\beta_{i}}\sum_{l=0}^{\beta_{i}-\beta_{j}-1}\left(
\frac{x_{j}}{x_{i}}\right)  ^{l}x^{\beta}-\kappa\sum_{\beta_{j}>\beta_{i}}%
\sum_{l=1}^{\beta_{j}-\beta_{i}}\left(  \frac{x_{i}}{x_{j}}\right)
^{l}x^{\beta}.
\]
The term $x^{\alpha}$ appears in the sum if (i) $\alpha=\beta$, (ii) (with
coefficient $\kappa$) for some $j$, $\beta_{i}>\beta_{j}$ and $\alpha
_{i}=\beta_{i}-l,\alpha_{j}=\beta_{j}+l$ with $0\leq l\leq\beta_{i}-\beta
_{j}-1$, (iii) (with coefficient $-\kappa$) for some $j$, $\beta_{i}<\beta
_{j}$ and $\alpha_{i}=\beta_{i}+l,\alpha_{j}=\beta_{j}-l$ with $1\leq
l\leq\beta_{j}-\beta_{i}$ (for (ii) and (iii) $\alpha_{k}=\beta_{k}$ for
$k\neq i,j$). In case (ii) $\alpha=\beta$ for $l=0$ and $\beta^{+}\succ
\alpha^{+}$ for $1\leq l\leq\beta_{i}-\beta_{j}-1$ by \cite[Lemma 8.2.3]{DX}.
In case (iii) $\alpha=\left(  i,j\right)  \beta$ for $l=\beta_{j}-\beta
_{i}=\alpha_{i}-\alpha_{j}>0$ and $\beta^{+}\succ\alpha^{+}$ for $1\leq
l\leq\beta_{i}-\beta_{j}-1$ as before.
\end{proof}

The Lemma will be used in analyzing the effect of $\mathcal{D}_{\ell\left(
\alpha\right)  }\left(  \kappa\right)  $ on $q_{\alpha}\left(  \kappa\right)
$, the polynomial defined in Theorem \ref{eigvect0}. The aim will be to show
it suffices to consider $\mathcal{D}_{\ell\left(  \alpha\right)  }\left(
\kappa\right)  \zeta_{\alpha}^{x}\left(  \kappa\right)  $. We point out that
for any given partition $\tau$ with $\tau_{2}\geq n$ there may be several
reasons why there can be no singular polynomial of isotype $\tau$, notably
there may be no eigenfunction of $\left\{  \mathcal{U}_{i}\left(  \kappa
_{0}\right)  :1\leq i\leq N\right\}  $ with the respective eigenvalues
$\left\{  \xi_{i}\left(  \lambda;\kappa_{0}\right)  :1\leq i\leq N\right\}  $,
where $\lambda$ is as specified in Theorem \ref{glambda}. Our proof singles
out one aspect, a certain nonvanishing coefficient of $\mathcal{D}%
_{\ell\left(  \lambda\right)  }\left(  \kappa_{0}\right)  g_{\lambda}$ which
applies to all cases.

\section{The two-part case}

In this section we consider the simplest case where $\tau=\left(  \tau
_{1},\tau_{2}\right)  $, $\kappa_{0}=-\frac{m}{n}$ with $\gcd\left(
m,n\right)  =1$ and $\tau_{2}=n$. By Corollary \ref{bigtau} $\tau_{1}=dn-1$
for some $d\geq2$ (since $\tau_{1}\geq\tau_{2}$). We will show that there is
no singular polynomial for $\kappa_{0}$ of isotype $\tau$. By Theorem
\ref{glambda}, if there exist singular polynomials for $\kappa_{0}$ of isotype
$\tau$ then there exists $g_{\lambda}=x^{\lambda}+\sum_{\beta\vartriangleleft
\lambda}A_{\beta}x^{\beta}$ with $\lambda=\left(  \left(  md\right)
^{n}\right)  ,\mathcal{D}_{i}\left(  \kappa_{0}\right)  g_{\lambda}=0$ and
$\mathcal{U}_{i}\left(  \kappa_{0}\right)  g_{\lambda}=\xi_{i}\left(
\lambda;\kappa_{0}\right)  g_{\lambda}$ for $1\leq i\leq N$. In fact,
$g_{\lambda}=\lim\limits_{\kappa\rightarrow\kappa_{0}}\zeta_{\lambda}%
^{x}\left(  \kappa\right)  $. This follows from Theorem \ref{eigvect0} because
there is no $\kappa_{0}$-critical pair $\left(  \lambda,\beta\right)  $ with
$\ell\left(  \beta\right)  \leq N$. (The background for this is detailed in
\cite{D2}; briefly $\mathrm{coef}\left(  \zeta_{\alpha}^{x},\beta\right)  $ is
independent of the number $M$ of variables provided $\max\left(  \ell\left(
\alpha\right)  ,\ell\left(  \beta\right)  \right)  \leq M$ (see equation
\ref{pimn}); also if $\ell\left(  \alpha\right)  \leq N<\ell\left(
\alpha\right)  +\left\vert \alpha\right\vert $ then not every factor of
$h\left(  \alpha,\kappa+1\right)  $ need appear as a pole of $\zeta_{\alpha
}^{x}$.) We start by computing $h\left(  \lambda,\kappa+1\right)  $ and
showing there is a unique $\beta$ so that $\left(  \lambda,\beta\right)  $ is
a $\kappa_{0}$-critical pair for $\lambda$ and $\ell\left(  \beta\right)
=N+1$.

For the rectangular diagram $\lambda=\left(  \left(  md\right)  ^{n}\right)  $
it is clear that $L\left(  \lambda;i,j\right)  =n-i$ for $1\leq i\leq n,1\leq
j\leq md$ so that $h\left(  \lambda,\kappa+1\right)  =\prod\nolimits_{i=1}%
^{n}\prod\nolimits_{j=1}^{md}\left(  \left(  n-i+1\right)  \kappa
+md+1-j\right)  =\prod\nolimits_{i=1}^{n}\prod\nolimits_{j=1}^{md}\left(
i\kappa+j\right)  $. Since $\gcd\left(  m,n\right)  =1$ the multiplicity of
$\left(  n\kappa+m\right)  $ in $h\left(  \lambda,\kappa+1\right)  $ is one,
occurring as $h\left(  \lambda,\kappa+1;1,md-m+1\right)  .$ The algorithm of
\cite{D3} yields $\beta=\left(  0^{n},m^{nd}\right)  $ for a $\kappa_{0}%
$-critical pair $\left(  \lambda,\beta\right)  $. Note $\ell\left(
\beta\right)  =n+nd=\tau_{2}+\left(  \tau_{1}+1\right)  =N+1$. Also recall the
easily proved rule: for any critical pair $\left(  \alpha,\gamma\right)  $ it
always holds that if $i>\ell\left(  \alpha\right)  $ and $\gamma_{i}=0$ then
$\gamma_{j}=0$ for all $j>i$; since $r\left(  \alpha,i\right)  =i=r\left(
\gamma,i\right)  $.

\begin{proposition}
\label{twoprt1}For $\lambda=\left(  \left(  md\right)  ^{n}\right)
,\kappa_{0}=-\frac{m}{n}$ with $\gcd\left(  m,n\right)  =1$ let $\beta=\left(
0^{n},m^{nd}\right)  $ then $\left(  \lambda,\beta\right)  $ is the unique
$\kappa_{0}$-critical pair for $\lambda$.
\end{proposition}

\begin{proof}
Suppose $\gamma\in\mathbb{N}_{0}^{M}$ for some $M\geq N$, and $\gamma$
satisfies the conditions $\lambda\trianglerighteq\gamma$ and $R_{i}:\left(
r\left(  \gamma,i\right)  -i\right)  m=\left(  \lambda_{i}-\gamma_{i}\right)
n$ for $1\leq i\leq M$ (as usual, define $\lambda_{i}=0$ for any
$i>\ell\left(  \lambda\right)  $, the equation is a restatement of $\left(
r\left(  \gamma,i\right)  -i\right)  \kappa_{0}+\left(  \lambda_{i}-\gamma
_{i}\right)  =0$). We must show $\gamma=\lambda$ or $\gamma=\beta$. Since
$\gcd\left(  m,n\right)  =1$ there exists $\eta\in\mathbb{N}_{0}^{M}$ so that
$\gamma=m\eta$ (componentwise; note $r\left(  \gamma,i\right)  =r\left(
\eta,i\right)  $ for each $i$). By condition $R_{n+1}$ we have $\left(
r\left(  \eta,n+1\right)  -n-1\right)  =-n\eta_{n+1}$ so that $\eta
_{n+1}=1-\frac{1}{n}\left(  r\left(  \eta,n+1\right)  -1\right)  \leq1$ and
thus $\eta_{n+1}=1$ or $\eta_{n+1}=0$.

If $\eta_{n+1}=1$ then $r\left(  \eta,n+1\right)  =1$, which implies $\eta
_{i}=0$ for $1\leq i\leq n$ and $\eta_{i}\leq1$ for $i>n+1.$ Since $\left\vert
\eta\right\vert =\frac{1}{m}\left\vert \gamma\right\vert =\frac{1}%
{m}\left\vert \lambda\right\vert =nd$ we see that $\eta^{+}=\left(
1^{nd}\right)  $ and in fact $\eta_{i}=1$ for $n+1\leq i\leq n\left(
d+1\right)  $, since $\eta_{j}=0$ and $\eta_{j+1}=1$ is impossible for $j>n$.
Thus $r\left(  1,\eta\right)  =nd+1$ and condition $R_{1}$ becomes $\left(
nd+1-1\right)  m=\left(  md-0\right)  n$. So $\gamma=\beta$; the other
conditions $R_{i}$ are verified similarly.

If $\eta_{n+1}=0$ then $r\left(  \eta,n+1\right)  =n+1$ and $\ell\left(
\eta\right)  =\ell\left(  \gamma\right)  =n$. But the conditions $\lambda
_{i}=\lambda_{1}$ for $1\leq i\leq n=\ell\left(  \lambda\right)
,\lambda\trianglerighteq\gamma$ and $\ell\left(  \gamma\right)  =\ell\left(
\lambda\right)  $ together imply $\gamma=\lambda$.
\end{proof}

\begin{corollary}
\label{twoprt2}The polynomial $\zeta_{\lambda}^{x}\left(  \kappa\right)  $ in
$N=nd+n-1$ variables has no pole at $\kappa_{0}$. The (hypothetical) singular
polynomial $g_{\lambda}=\lim\limits_{\kappa\rightarrow\kappa_{0}}%
\zeta_{\lambda}^{x}\left(  \kappa\right)  $.
\end{corollary}

\begin{proof}
By Theorem \ref{eigvect0} $g_{\lambda}=\lim\limits_{\kappa\rightarrow
\kappa_{0}}\zeta_{\lambda}^{x}\left(  \kappa\right)  $, since there is no
$\gamma\in\mathbb{N}_{0}^{N}$ so that $\left(  \lambda,\gamma\right)  $ is
$\kappa_{0}$-critical. By \cite[Theorem 4.8]{D2} $\zeta_{\lambda}^{x}\left(
\kappa\right)  $ in $nd+n-1$ variables has no pole at $\kappa_{0}$ (since
$N=nd+n-1<\ell\left(  \beta\right)  $).
\end{proof}

This implies $\mathcal{D}_{n}\left(  \kappa_{0}\right)  g_{\lambda}%
=\lim\limits_{\kappa\rightarrow\kappa_{0}}\mathcal{D}_{n}\left(
\kappa\right)  \zeta_{\lambda}^{x}\left(  \kappa\right)  $. In the notation of
Section 3 (noting $N+1-r\left(  \lambda,n\right)  =\left(  nd+n-1\right)
+1-n=nd$)
\[
\mathcal{D}_{n}\left(  \kappa\right)  \zeta_{\lambda}^{x}\left(
\kappa\right)  =\frac{md}{\kappa+md}d\left(  n\kappa+m\right)  \theta_{n}%
^{-1}\zeta_{\widetilde{\lambda}}^{x}\left(  \kappa\right)  ,
\]
where $\widetilde{\lambda}=\left(  md-1,\left(  md\right)  ^{n-1}\right)  $.
We will show that the coefficient of $x^{\gamma}$ in the equation does not
converge to zero as $\kappa\rightarrow\kappa_{0}$, where $\gamma=\theta
_{n}\left(  m-1,0^{n-1},m^{nd-1}\right)  =\left(  0^{n-1},m-1,m^{nd-1}\right)
$ and $\ell\left(  \gamma\right)  =N$. The following is similar to Proposition
\ref{twoprt1}.

\begin{proposition}
\label{uniqb}For $\widetilde{\lambda}=\left(  md-1,\left(  md\right)
^{n-1}\right)  ,\kappa_{0}=-\frac{m}{n}$ with $\gcd\left(  m,n\right)  =1$ let
$\beta=\left(  m-1,0^{n-1},m^{nd-1}\right)  $ then $\left(  \widetilde
{\lambda},\beta\right)  $ is the unique $\kappa_{0}$-critical pair for
$\widetilde{\lambda}$.
\end{proposition}

\begin{proof}
Suppose $\gamma\in\mathbb{N}_{0}^{M}$ for some $M\geq N$, and $\gamma$
satisfies the conditions $\widetilde{\lambda}\trianglerighteq\gamma$ and
$R_{i}:\left(  r\left(  \gamma,i\right)  -r\left(  \widetilde{\lambda
},i\right)  \right)  m=\left(  \widetilde{\lambda}_{i}-\gamma_{i}\right)  n$
for $1\leq i\leq M$. The conditions $R_{i}$ specialize to:%
\begin{align*}
\left(  r\left(  \gamma,1\right)  -n\right)  m  & =\left(  md-1-\gamma
_{1}\right)  n,\text{ for }i=1,\\
\left(  r\left(  \gamma,i\right)  -i+1\right)  m  & =\left(  md-\gamma
_{i}\right)  n,\text{ for }2\leq i\leq n,\\
\left(  r\left(  \gamma,i\right)  -i+1\right)  m  & =-\gamma_{i}n,\text{ for
}i>n.
\end{align*}
Thus $\gamma_{1}\equiv m-1\operatorname{mod}m$ and $\gamma_{i}\equiv
0\operatorname{mod}m$ for $i\geq2$. Consider the condition $R_{n+1}:\left(
r\left(  \eta,n+1\right)  -n-1\right)  m=-n\gamma_{n+1}$, equivalent to
$\gamma_{n+1}/m=1-\left(  r\left(  \gamma,n+1\right)  -1\right)  /n\leq1$.
Thus $\gamma_{n+1}=m$ or $\gamma_{n+1}=0$. If $\gamma_{n+1}=m$ then $r\left(
\gamma,n+1\right)  =1$, implying that $\gamma_{i}<m$ for $1\leq i\leq n$. The
congruence conditions imply $\gamma_{1}=m-1$ and $\gamma_{i}=0$ for $2\leq
i\leq n$. Arguing similarly to Proposition \ref{twoprt1} we see that
$\gamma_{i}=m$ or $\gamma_{i}=0$ for $i\geq n+1$; the condition $\left\vert
\gamma\right\vert =nmd-1$ shows that $\gamma=\left(  m-1,0^{n-1}%
,m^{nd-1}\right)  $. To verify $R_{1}$ note $r\left(  \gamma,1\right)  =nd$ so
$\left(  nd-n\right)  m=\left(  \left(  md-1\right)  -\left(  m-1\right)
\right)  n=mn\left(  d-1\right)  $ is satisfied. For $2\leq i\leq n$ we have
$r\left(  \gamma,i\right)  =nd+i-1$ and $\gamma_{i}=0$ so the condition
$R_{i}$, namely, $\left(  nd+i-1-i+1\right)  m=mdn$ is satisfied.

If $\gamma_{n+1}=0$ then $\ell\left(  \gamma\right)  =n$; for this particular
$\widetilde{\lambda}$ the conditions $\left(  \widetilde{\lambda}\right)
^{+}\succeq\gamma^{+}$ and $\ell\left(  \gamma\right)  =n$ together imply
$\left(  \widetilde{\lambda}\right)  ^{+}=\gamma^{+}$; then $\widetilde
{\lambda}\succeq\gamma$ implies $\widetilde{\lambda}=\gamma$.
\end{proof}

By Lemma \ref{hookd} $h\left(  \widetilde{\lambda},\kappa+1\right)  =h\left(
\lambda,\kappa+1\right)  /\left(  \kappa+md\right)  $ so $\left(
n\kappa+m\right)  $ has multiplicity one in $h\left(  \widetilde{\lambda
},\kappa+1\right)  $. Next we will show $\mathrm{coef}\left(  \zeta
_{\widetilde{\lambda}}\left(  \kappa\right)  ,\beta\right)  $ has a simple
pole at $\kappa_{0}$, where $\beta$ is defined in the Proposition.

For $w\in S_{N}$ and $\alpha\in\mathbb{N}_{0}^{N}$ since $w\left(  x^{\alpha
}\right)  =x^{w\alpha}$ the transformation property $\mathrm{coef}\left(
p,\alpha\right)  =\mathrm{coef}\left(  wp,w\alpha\right)  $ holds for any
polynomial $p$.

\begin{theorem}
Suppose $\kappa_{0}=-\frac{m}{n}$ with $\gcd\left(  m,n\right)  =1$ and
$\tau=\left(  dn-1,n\right)  $ with $d\geq2,n\geq2$ so that $N=\left(
d+1\right)  n-1$, then there are no singular polynomials for $\kappa_{0}$ of
isotype $\tau$.
\end{theorem}

\begin{proof}
By Corollary \ref{twoprt2} if there is a singular polynomial of isotype $\tau$
for the singular value $\kappa_{0}$ then $g_{\lambda}=\lim\limits_{\kappa
\rightarrow\kappa_{0}}\zeta_{\lambda}^{x}\left(  \kappa\right)  $ is singular,
where $\lambda=\left(  \left(  md\right)  ^{n}\right)  $. By Proposition
\ref{zdiff} $\mathcal{D}_{n}\left(  \kappa\right)  \zeta_{\lambda}^{x}\left(
\kappa\right)  =\frac{md}{\kappa+md}d\left(  n\kappa+m\right)  \theta_{n}%
^{-1}\zeta_{\widetilde{\lambda}}^{x}\left(  \kappa\right)  $. Let
$\beta=\left(  m-1,0^{n-1},m^{nd-1}\right)  $ so that $\left(  \widetilde
{\lambda},\beta\right)  $ is the unique $\kappa_{0}$-critical pair for
$\widetilde{\lambda}$. It is crucial that $\ell\left(  \beta\right)
=n+nd-1=N$. By Lemma \ref{kpole} $\mathrm{coef}\left(  \zeta_{\widetilde
{\lambda}}^{x}\left(  \kappa\right)  ,\beta\right)  =\frac{f\left(
\kappa\right)  }{n\kappa+m}$ where $f\left(  \kappa\right)  \in\mathbb{Q}%
\left(  \kappa\right)  $, $f\left(  \kappa\right)  $ has no pole at
$\kappa_{0}=-\frac{m}{n}$ and $f\left(  \kappa_{0}\right)  \neq0$. Note
$\theta_{n}^{-1}\beta=\left(  0^{n-1},m-1,m^{nd-1}\right)  $. Thus%
\[
\mathrm{coef}\left(  \mathcal{D}_{n}\left(  \kappa\right)  \zeta_{\lambda}%
^{x}\left(  \kappa\right)  ,\theta_{n}^{-1}\beta\right)  =\dfrac{md^{2}\left(
n\kappa+m\right)  f\left(  \kappa\right)  }{\left(  \kappa+md\right)  \left(
n\kappa+m\right)  }=\dfrac{md^{2}f\left(  \kappa\right)  }{\kappa+md}%
\]
and%
\[
\mathrm{coef}\left(  \mathcal{D}_{n}\left(  \kappa_{0}\right)  g_{\lambda
},\theta_{n}^{-1}\beta\right)  =\lim\limits_{\kappa\rightarrow\kappa_{0}%
}\dfrac{md^{2}f\left(  \kappa\right)  }{\kappa+md}\neq0,
\]
and so $g_{\lambda}$ is not singular for $\kappa_{0}$, a contradiction.
\end{proof}

This argument will serve as the key ingredient for the general case $\tau$.

\section{The general case}

In this section we consider the singular value $\kappa_{0}=-\frac{m}{n}$ with
$\gcd\left(  m,n\right)  =1$ and $2\leq n\leq N$ for the isotype $\tau$, where
$\tau$ has two or more parts and $\tau_{2}>n$. Let $l=\ell\left(  \tau\right)
\geq2$. We assume there exists a singular polynomial for $\kappa_{0}$ of
isotype $\tau$ and will eventually arrive at a contradiction. By Corollary
\ref{bigtau} there are integers $d_{1},d_{2},\ldots,d_{l-1}$ so that $\tau
_{i}=d_{i}n-1$ for $1\leq i<l$. Because $\tau$ is a partition it follows that
$d_{1}\geq d_{2}\geq\ldots\geq d_{l-1}$. By hypothesis $\tau_{1}\geq\tau
_{2}>n$ so that $d_{1}\geq2$, and $d_{2}\geq2$ if $l\geq3$. By Theorem
\ref{glambda} there is a corresponding partition $\lambda$ and a singular
polynomial $g_{\lambda}=x^{\lambda}+\sum_{\beta\vartriangleleft\lambda
}A_{\beta}x^{\beta}$ with $\mathcal{U}_{i}\left(  \kappa_{0}\right)
g_{\lambda}=\xi_{i}\left(  \lambda;\kappa_{0}\right)  g_{\lambda}$ for $1\leq
i\leq N$. The computations are expressed in terms of (with $1\leq i\leq l-1$):%
\begin{align*}
t_{i}  &  :=\sum_{j=1}^{l-i}d_{j}=\sum_{j=1}^{l-i}\frac{\tau_{j}+1}{n},\\
p_{i}  &  :=\sum_{j=l+1-i}^{l}\tau_{j},\\
\lambda &  :=\left(  \left(  mt_{1}\right)  ^{\tau_{l}},\left(  mt_{2}\right)
^{\tau_{l-1}},\ldots,\left(  mt_{l-1}\right)  ^{\tau_{2}},0^{\tau_{1}}\right)
,\\
\gamma &  :=\left(  \left(  mt_{1}\right)  ^{\tau_{l}},\left(  mt_{2}\right)
^{\tau_{l-1}},\ldots,\left(  mt_{l-1}\right)  ^{\tau_{2}-n},0^{n-1}%
,m-1,m^{\tau_{1}}\right)  ,\\
\alpha &  :=\left(  \left(  mt_{1}\right)  ^{\tau_{l}},\left(  mt_{2}\right)
^{\tau_{l-1}},\ldots,\left(  mt_{l-1}\right)  ^{\tau_{2}-n},0^{n-1}%
,m,m^{\tau_{1}}\right)  .
\end{align*}
Also let $p_{0}=0,p_{l}=N$. Note that $p_{l-1}=\ell\left(  \lambda\right)
=N-\tau_{1}$, $t_{l-1}=d_{1}\geq2$ and $\left\vert \gamma\right\vert
+1=\left\vert \alpha\right\vert =\left\vert \lambda\right\vert $ because
$\left(  \tau_{1}+1\right)  m=nd_{1}m=nmt_{l-1}$. By Theorem \ref{eigvect0}
there exists a polynomial $q_{\lambda}\left(  \kappa\right)  =\zeta_{\lambda
}^{x}\left(  \kappa\right)  +\sum\left\{  B_{\beta}\left(  \kappa\right)
\zeta_{\beta}^{x}\left(  \kappa\right)  :\beta\in C\left(  \lambda,\kappa
_{0}\right)  \right\}  $ which has no pole at $\kappa_{0}$ and $\lim
\limits_{\kappa\rightarrow\kappa_{0}}q_{\lambda}\left(  \kappa\right)
=g_{\lambda}$; the coefficients $B_{\beta}\left(  \kappa\right)  \in
\mathbb{Q}\left(  \kappa,v\right)  $ and $C\left(  \lambda,\kappa_{0}\right)
$ is the set of $\beta$ such that $\left(  \lambda,\beta\right)  $ is a
$\kappa_{0}$-critical pair.

We will show $\lim\limits_{\kappa\rightarrow\kappa_{0}}\mathcal{D}%
_{\ell\left(  \lambda\right)  }\left(  \kappa\right)  q_{\lambda}\left(
\kappa\right)  \neq0$ by showing $\lim\limits_{\kappa\rightarrow\kappa_{0}%
}\mathrm{coef}\left(  \mathcal{D}_{\ell\left(  \lambda\right)  }\left(
\kappa\right)  q_{\lambda}\left(  \kappa\right)  ,\gamma\right)  \neq0$. The
proof has two parts: firstly we show that $\mathrm{coef}\left(  \mathcal{D}%
_{\ell\left(  \lambda\right)  }\left(  \kappa\right)  q_{\lambda}\left(
\kappa\right)  ,\gamma\right)  =\mathrm{coef}\left(  \mathcal{D}_{\ell\left(
\lambda\right)  }\left(  \kappa\right)  \zeta_{\lambda}^{x}\left(
\kappa\right)  ,\gamma\right)  $ and secondly we use the Insertion Theorem
\ref{coefs1} and the result from the previous section.

\begin{lemma}
\label{dxg}Suppose $\delta\in\mathbb{N}_{0}^{N},\lambda\trianglerighteq\delta$
and $\mathrm{coef}\left(  \mathcal{D}_{\ell\left(  \lambda\right)  }\left(
\kappa\right)  x^{\delta},\gamma\right)  \neq0$ then $\delta\trianglerighteq
\alpha$.
\end{lemma}

\begin{proof}
By construction $\mathrm{coef}\left(  \mathcal{D}_{\ell\left(  \lambda\right)
}\left(  \kappa\right)  x^{\delta},\gamma\right)  =\mathrm{coef}\left(
x_{\ell\left(  \lambda\right)  }\mathcal{D}_{\ell\left(  \lambda\right)
}\left(  \kappa\right)  x^{\delta},\alpha\right)  $. By Lemma \ref{dcoefs}
$\delta=\alpha$ or $\delta^{+}\succ\alpha^{+}$ or $\delta=\left(  \ell\left(
\lambda\right)  ,j\right)  \alpha$ with $\alpha_{\ell\left(  \lambda\right)
}>\alpha_{j}$; but in the latter case $j<\ell\left(  \lambda\right)  $ (in
fact $\ell\left(  \lambda\right)  -n\leq j<\ell\left(  \lambda\right)  $ and
$\alpha_{j}=0$) so that $\delta\succ\alpha$. Thus $\delta\trianglerighteq
\alpha$.
\end{proof}

To complete the first part of the argument we need only show that there is no
$\kappa_{0}$-critical pair $\left(  \lambda,\beta\right)  $ with
$\beta\trianglerighteq\alpha$.

\begin{theorem}
Suppose $\beta\in\mathbb{N}_{0}^{M}$ (with some $M\geq N$), $\lambda
\trianglerighteq\beta\trianglerighteq\alpha$ and $\beta$ satisfies the rank
equation $\left(  r\left(  \beta,i\right)  -i\right)  m=\left(  \lambda
_{i}-\beta_{i}\right)  n$ for $i\geq1$ then $\beta=\lambda$.
\end{theorem}

\begin{proof}
The rank equation and definition of $\alpha$ imply $m|\beta_{i}$ and
$m|\alpha_{i}$ for each $i$ ; since the definition of $\vartriangleright$
implies that $\mu\vartriangleright\nu$ if and only if $m\mu\vartriangleright
m\nu$ for arbitrary compositions $\mu,\nu$ (where $\left(  m\mu\right)
_{i}:=m\mu_{i}$) we will assume that $m=1$ in the rest of the proof. Since
$\beta$ is trapped between $\lambda=\left(  t_{1}^{\tau_{l}},t_{2}^{\tau
_{l-1}},\ldots,t_{l-1}^{\tau_{2}},0^{\tau_{1}}\right)  $ and $\alpha=\left(
t_{1}^{\tau_{l}},t_{2}^{\tau_{l-1}},\ldots,t_{l-1}^{\tau_{2}-n},0^{n-1}%
,1^{\tau_{1}+1}\right)  $ we deduce that $\beta^{+}$ agrees with $\lambda$ in
the first $N-\tau_{1}-n$ entries, that is, $\left(  \beta^{+}\right)
_{i}=\lambda_{i}=\alpha_{i}$ for $1\leq i\leq N-\tau_{1}-n=p_{l-1}-n$. None of
the entries of $\beta$ equal to some $t_{j}$ can \textquotedblleft move to the
left\textquotedblright\ (lower index). For $j,k$ with $1\leq j\leq k\leq l-1$
suppose that the first appearance (least index) of $t_{k}$ in $\beta$ is at an
index $i$ with $\lambda_{i}=t_{j}$, that is, $p_{j-1}+1\leq i\leq p_{j}$ ,
then $r\left(  \beta,i\right)  =p_{k-1}+1$ and the rank equation implies%
\begin{align*}
i  &  =r\left(  \beta,i\right)  -n(\lambda_{i}-\beta_{i})=p_{k-1}+1-n\left(
t_{j}-t_{k}\right) \\
&  =p_{k-1}+1-n\sum_{s=l-k+1}^{l-j}d_{s}=p_{k-1}+1-\sum_{s=l-k+1}^{l-j}\left(
\tau_{s}+1\right)  .
\end{align*}
Furthermore
\begin{align*}
0  &  \leq i-\left(  p_{j-1}+1\right)  =p_{k-1}-p_{j-1}-\sum_{s=l-k+1}%
^{l-j}\left(  \tau_{s}+1\right) \\
&  =\sum_{s=l+2-k}^{l+1-j}\tau_{s}-\sum_{s=l-k+1}^{l-j}\tau_{s}-\left(
k-j\right)  =\tau_{l+1-j}-\tau_{l-k+1}-\left(  k-j\right) \\
&  \leq j-k\leq0.
\end{align*}
The inequality $\tau_{l+1-j}-\tau_{l-k+1}\leq0$ holds because $\tau$ is a
partition and $l-k+1\leq l-j+1$ by hypothesis. The chain of inequalities shows
that $j=k$ and $i=p_{k-1}+1$ (the possibility $i>p_{k-1}+1$ has not yet been
ruled out).

The key to the argument is the value of $\beta_{\ell\left(  \lambda\right)
+1}$ (recall $\ell\left(  \lambda\right)  =p_{l-1}$). The case $\beta
_{p_{l-1}+1}=t_{j}$ is impossible for $1\leq j\leq l-1$; indeed suppose
$\beta_{p_{l-1}+1}=t_{j}$ then $r\left(  \beta,p_{l-1}+1\right)  \geq
p_{j-1}+1$ and the rank equation is $p_{l-1}+1=r\left(  \beta,p_{l-1}%
+1\right)  +nt_{j}$ (note $\lambda_{p_{l-1}+1}=0$) thus
\begin{align*}
0  &  \leq\left(  p_{l-1}+1-nt_{j}\right)  -\left(  p_{j-1}+1\right) \\
&  =\tau_{l+1-j}-\tau_{1}-\left(  l-j\right)  \leq j-l<0.
\end{align*}
which is a contradiction (the calculation is similar to the previous one,
replacing $k$ by $l$ and $t_{k}$ by $t_{l}=0$). The condition $\beta
^{+}\preceq\lambda$ now implies that $\beta_{p_{l-1}+1}<t_{l-1}$ and
$\#\left\{  j:\beta_{j}>\beta_{p_{l-1}+1}\right\}  \geq\ell\left(
\lambda\right)  -n=p_{l-1}-n$, hence $r\left(  \beta,p_{l-1}+1\right)  \geq
p_{l-1}-n+1$. The rank equation is $-n\beta_{p_{l-1}+1}=r\left(  \beta
,p_{l-1}+1\right)  -\left(  p_{l-1}+1\right)  \geq-n$ and so $\beta
_{p_{l-1}+1}\leq1$.

Suppose $\beta_{p_{l-1}+1}=0$ then $r\left(  \beta,p_{l-1}+1\right)  =\left(
p_{l-1}+1\right)  $ which implies $\beta_{i}=0$ for $i>p_{l-1}+1$ and
$\beta_{i}>0$ for $i\leq p_{l-1}$. Since $\beta^{+}$ differs from $\lambda$ in
at most the last $n$ entries and $\lambda_{i}=t_{l-1}$ for $\ell\left(
\lambda\right)  -n<i\leq\ell\left(  \lambda\right)  $ it follows that
$\beta^{+}=\lambda$ (note if $\left(  t_{l-1}^{n}\right)  \succeq\mu$ where
$\mu$ is a partition and $\ell\left(  \mu\right)  =n$ then $\mu=\left(
t_{l-1}^{n}\right)  $). Since the entries of $\beta$ can not move to the left,
$\beta=\lambda$; in detail, argue inductively that the only possible value for
$\beta_{i}$ when $1\leq i\leq p_{1}$ is $t_{1}$, then the only possible value
when $p_{1}+1\leq i\leq p_{2}$ is $t_{2}$, and so on. (If $l=2$ then this
argument is not needed.)

Suppose $\beta_{p_{l-1}+1}=1$. In this part replace the bound $\beta
\trianglerighteq\alpha$ by $\beta\trianglerighteq\alpha^{\prime}:=\left(
t_{1}^{\tau_{l}},t_{2}^{\tau_{l-1}},\ldots,t_{l-1}^{\tau_{2}-n},0^{n}%
,1^{\tau_{1}+1}\right)  $, a weaker restriction since $\alpha\vartriangleright
\alpha^{\prime}$ (note that $\left(  \lambda,\alpha^{\prime}\right)  $ is a
$\left(  -\frac{1}{n}\right)  $-critical pair). We will show that
$\beta=\alpha^{\prime}$, which contradicts the assumption $\ell\left(
\beta\right)  =N$. Recall $t_{j}>t_{l-1}\geq2$ for $1\leq j<l-1$, by the
hypothesis $\tau_{1}\geq\tau_{2}>n$. The rank equation yields $r\left(
\beta,p_{l-1}+1\right)  =\left(  p_{l-1}+1\right)  -n\beta_{p_{l-1}+1}%
=p_{l-1}+1-n$. For $i<p_{l-1}+1$ this implies $\beta_{i}=t_{j}$ for some $j$
or $\beta_{i}<1$, that is, $\beta_{i}=0$; for $i>p_{l-1}+1$ the rank implies
$\beta_{i}=t_{j}$ for some $j$ or $\beta_{i}\leq1$. This forces the values of
$\beta$ other than $\left(  t_{1}^{\tau_{l}},t_{2}^{\tau_{l-1}},\ldots
,t_{l-1}^{\tau_{2}-n}\right)  $ to be $0$ or $1$, that is, $\left(  \beta
^{+}\right)  _{i}\leq1$ for $i>p_{l-1}-n$. The condition $\left\vert
\lambda\right\vert =\left\vert \beta\right\vert $ shows that $\#\left\{
j:\beta_{j}=1\right\}  =nt_{l-1}=\tau_{1}+1$. Since $\beta_{i}=1$ is ruled out
for $i\leq p_{l-1}$ it follows that $\beta_{i}=1$ for $p_{l-1}+1\leq i\leq
p_{l-1}+\tau_{1}+1=N+1$; indeed suppose the $j^{th}$ occurrence of $1$ in
$\beta$ is at index $i$, that is, $r\left(  \beta,i\right)  =p_{l-1}+j-n$,
$\beta_{i}=1$ and $i>p_{l-1}$ then the rank equation implies $r\left(
\beta,i\right)  -i=\left(  p_{l-1}+j-n\right)  -i=-n\beta_{i}=-n$ thus
$i=p_{l-1}+j$, for $1\leq j\leq\tau_{1}+1$. The $n$ remaining values of
$\beta_{i}$ (for $1\leq i\leq N+1$) are zero, and so $\beta^{+}=\left(
\alpha^{\prime}\right)  ^{+}.$ The condition $\beta\trianglerighteq
\alpha^{\prime}$ implies $\beta\succeq\alpha^{\prime}$ (by definition) which
shows $\beta_{i}=\alpha_{i}^{\prime}$ for $1\leq i\leq p_{l-1}-n$ (if
$p_{j-1}+1\leq i\leq p_{j}$ and $j<l-1$ then $\beta_{i}=t_{j}$ and if
$p_{l-2}+1\leq i\leq p_{l-1}-n$ then $\beta_{i}=t_{l-1}$). Thus $\beta_{i}=0$
for $p_{l-1}-n<i\leq p_{l-1}$ and $\beta=\alpha^{\prime}$. The proof is
finished since $\ell\left(  \alpha^{\prime}\right)  =N+1.$
\end{proof}

\begin{corollary}
$\mathrm{coef}\left(  \mathcal{D}_{\ell\left(  \lambda\right)  }\left(
\kappa\right)  q_{\lambda}\left(  \kappa\right)  ,\gamma\right)
=\mathrm{coef}\left(  \mathcal{D}_{\ell\left(  \lambda\right)  }\left(
\kappa\right)  \zeta_{\lambda}^{x}\left(  \kappa\right)  ,\gamma\right)  $.
\end{corollary}

\begin{proof}
Suppose $\beta\in C\left(  \lambda,\kappa_{0}\right)  $. If $\mathrm{coef}%
\left(  \mathcal{D}_{\ell\left(  \lambda\right)  }\left(  \kappa\right)
\zeta_{\beta}^{x}\left(  \kappa\right)  ,\gamma\right)  \neq0$ then by Lemma
\ref{dxg} there exists $\delta\in\mathbb{N}_{0}^{N}$ such that $\beta
\trianglerighteq\delta\trianglerighteq\alpha$, which contradicts the Theorem.
Hence $\mathrm{coef}\left(  \mathcal{D}_{\ell\left(  \lambda\right)  }\left(
\kappa\right)  \zeta_{\beta}^{x}\left(  \kappa\right)  ,\gamma\right)  =0$ for
each $\beta\in C\left(  \lambda,\kappa_{0}\right)  $.
\end{proof}

\begin{example}
In the context of the Theorem there may well be compositions $\beta$ other
than $\alpha^{\prime}$ for which $\left(  \lambda,\beta\right)  $ is a
$\kappa_{0}$-critical pair. Suppose $N=10,~\tau=\left(  3,3,3,1\right)  $ and
$\kappa_{0}=-\frac{1}{2}$, then $\lambda=\left(  6,4^{3},2^{3}\right)  $ and
$\alpha^{\prime}=\left(  6,4^{3},2,0,0,1^{4}\right)  $; the multiplicity of
$\left(  2\kappa+1\right)  $ in $h\left(  \lambda,\kappa+1\right)  $ is $3$.
Among other compositions $\beta$ with $\left(  \lambda,\beta\right)  $ being
$\left(  -\frac{1}{2}\right)  $-critical are $\left(  6,1^{3},2^{3}%
,3^{3}\right)  $ and $\left(  6,0^{3},2,4^{2},1^{4},4\right)  $; the latter is
a permutation of $\alpha^{\prime}$. For another example take $N=14,~\tau
=\left(  8,6\right)  $ and $\kappa_{0}=-\frac{1}{3}$, then $\lambda=\left(
3^{6}\right)  $ and $\alpha^{\prime}=\left(  3^{3},0^{3},1^{9}\right)  $; the
multiplicity of $\left(  3\kappa+1\right)  $ in $h\left(  \lambda
,\kappa+1\right)  $ is $2$ and both $\left(  1^{6},2^{6}\right)  $ and
$\left(  1^{3},0^{3},2^{6},1^{3}\right)  $ form $\left(  -\frac{1}{3}\right)
$-critical pairs with $\lambda$. The algorithm of \cite{D3} was used to
produce the $\beta$'s.
\end{example}

Let $k=\ell\left(  \lambda\right)  =N-\tau_{1}$ and from Definition
\ref{cyclic} let $\theta_{k}=\left(  1,2\right)  \ldots\left(  k-1,k\right)
\in S_{N}$, a cyclic shift and let
\begin{align*}
\widetilde{\lambda}  &  =\left(  mt_{l-1}-1,\left(  mt_{1}\right)  ^{\tau_{l}%
},\left(  mt_{2}\right)  ^{\tau_{l-1}},\ldots,\left(  mt_{l-1}\right)
^{\tau_{2}-1},0^{\tau_{1}}\right)  ,\\
\widetilde{\alpha}  &  =\left(  m-1,\left(  mt_{1}\right)  ^{\tau_{l}},\left(
mt_{2}\right)  ^{\tau_{l-1}},\ldots,\left(  mt_{l-1}\right)  ^{\tau_{2}%
-n},0^{n-1},m^{\tau_{1}}\right)  ,
\end{align*}
so that $\widetilde{\alpha}=\theta_{k}\gamma$. By Proposition \ref{zdiff}
\[
\mathcal{D}_{\ell\left(  \lambda\right)  }\left(  \kappa\right)
\zeta_{\lambda}^{x}\left(  \kappa\right)  =\dfrac{mt_{l-1}}{\kappa+mt_{l-1}%
}\left(  \left(  N+1-k\right)  \kappa+mt_{l-1}\right)  \theta_{k}^{-1}%
\zeta_{\widetilde{\lambda}}^{x}\left(  \kappa\right)
\]
and $\left(  N+1-k\right)  \kappa+mt_{l-1}=\left(  \tau_{1}+1\right)
\kappa+md_{1}=\left(  n\kappa+m\right)  d_{1}$ (recall $\tau_{1}+1=nd_{1}$).
Thus%
\begin{align*}
\mathrm{coef}\left(  \mathcal{D}_{\ell\left(  \lambda\right)  }\left(
\kappa\right)  \zeta_{\lambda}^{x}\left(  \kappa\right)  ,\gamma\right)   &
=\dfrac{md_{1}^{2}}{\kappa+md_{1}}\left(  n\kappa+m\right)  \mathrm{coef}%
\left(  \theta_{k}^{-1}\zeta_{\widetilde{\lambda}}^{x}\left(  \kappa\right)
,\gamma\right) \\
&  =\dfrac{md_{1}^{2}}{\kappa+md_{1}}\left(  n\kappa+m\right)  \mathrm{coef}%
\left(  \zeta_{\widetilde{\lambda}}^{x}\left(  \kappa\right)  ,\theta
_{k}\gamma\right)  .
\end{align*}
We finish the argument by using the Insertion Theorem \ref{coefs1}. Let%
\begin{align*}
\mu &  =\left(  \left(  mt_{1}\right)  ^{\tau_{l}},\left(  mt_{2}\right)
^{\tau_{l-1}},\ldots,\left(  mt_{l-1}\right)  ^{\tau_{2}-n}\right)  ,\\
\nu &  =\left(  mt_{l-1}-1,\left(  mt_{l-1}\right)  ^{n-1},0^{\tau_{1}%
}\right)  ,\\
\sigma &  =\left(  m-1,0^{n-1},m^{\tau_{1}}\right)  ,
\end{align*}
then $\widetilde{\lambda}=\iota\left(  1,\mu\right)  \nu$ and $\widetilde
{\alpha}=\iota\left(  1,\mu\right)  \sigma$.

\begin{lemma}
$\lim\limits_{\kappa\rightarrow\kappa_{0}}\mathrm{coef}\left(  \mathcal{D}%
_{\ell\left(  \lambda\right)  }\left(  \kappa\right)  \zeta_{\lambda}%
^{x}\left(  \kappa\right)  ,\gamma\right)  \neq0$.
\end{lemma}

\begin{proof}
Let $p=mt_{l-1}=md_{1}$, and $M=n+\tau_{1}$ then $\nu,\sigma\in I_{1,p}%
^{\left(  M\right)  }$. By Theorem \ref{coefs1} $\mathrm{coef}\left(
\zeta_{\widetilde{\lambda}}^{x}\left(  \kappa\right)  ,\widetilde{\alpha
}\right)  =\mathrm{coef}\left(  \zeta_{\nu}^{x}\left(  \kappa\right)
,\sigma\right)  $ and by Proposition \ref{uniqb} and Lemma \ref{kpole}
$\mathrm{coef}\left(  \zeta_{\nu}^{x}\left(  \kappa\right)  ,\sigma\right)  $
has a simple pole at $\kappa=\kappa_{0}$ (this is the same argument used in
the previous section). Thus there exists $f\left(  \kappa\right)
\in\mathbb{Q}\left(  \kappa\right)  $ so that $\mathrm{coef}\left(  \zeta
_{\nu}^{x}\left(  \kappa\right)  ,\sigma\right)  =\frac{f\left(
\kappa\right)  }{n\kappa+m}$ and $f\left(  \kappa_{0}\right)  \neq0$. To
conclude,%
\begin{align*}
\mathrm{coef}\left(  \mathcal{D}_{\ell\left(  \lambda\right)  }\left(
\kappa\right)  \zeta_{\lambda}^{x}\left(  \kappa\right)  ,\gamma\right)   &
=\dfrac{md_{1}^{2}}{\kappa+md_{1}}\left(  n\kappa+m\right)  \mathrm{coef}%
\left(  \zeta_{\widetilde{\lambda}}^{x}\left(  \kappa\right)  ,\widetilde
{\alpha}\right) \\
&  =\dfrac{md_{1}^{2}}{\kappa+md_{1}}f\left(  \kappa\right)
\end{align*}
which has a nonzero limit at $\kappa_{0}.$
\end{proof}

\begin{theorem}
Suppose $\kappa_{0}=-\frac{m}{n}$ with $\gcd\left(  m,n\right)  =1$ and $\tau$
is a partition of $N$ such that $n|\left(  \tau_{i}+1\right)  $ for $1\leq
i<\ell\left(  \tau\right)  $. If $\tau_{2}\geq n$ then there are no singular
polynomials for $\kappa_{0}$ of isotype $\tau$.
\end{theorem}

\begin{proof}
For $\tau_{2}>n$ if there is a singular polynomial of isotype $\tau$ for the
singular value then $\lim\limits_{\kappa\rightarrow\kappa_{0}}\mathcal{D}%
_{\ell\left(  \lambda\right)  }\left(  \kappa\right)  q_{\lambda}\left(
\kappa\right)  =0$ for the polynomial $q_{\lambda}\left(  \kappa\right)  $
described above. But $\lim\limits_{\kappa\rightarrow\kappa_{0}}\mathrm{coef}%
\left(  \mathcal{D}_{\ell\left(  \lambda\right)  }\left(  \kappa\right)
q_{\lambda}\left(  \kappa\right)  ,\gamma\right)  =\lim\limits_{\kappa
\rightarrow\kappa_{0}}\mathrm{coef}\left(  \mathcal{D}_{\ell\left(
\lambda\right)  }\left(  \kappa\right)  \zeta_{\lambda}^{x}\left(
\kappa\right)  ,\gamma\right)  \neq0$, and so these singular polynomials do
not exist. The case $\tau_{2}=n,\ell\left(  \tau\right)  =2$ was done in the
previous section.
\end{proof}

\section{Concluding remarks}

Together with the results of \cite[Theorem 2.7]{D2} we have a complete
description of singular polynomials for the group $S_{N}$. For each pair
$\left(  m_{0},n_{0}\right)  \in\mathbb{N}^{2}$ with $2\leq n_{0}\leq N$ and
$\frac{m_{0}}{n_{0}}\notin\mathbb{N}$, let $d=\gcd\left(  m_{0},n_{0}\right)
,m=\frac{m_{0}}{d},n=\frac{n_{0}}{d}$, then there is a unique irreducible
$S_{N}$-module of singular polynomials for the singular value $\kappa
_{0}=-\frac{m}{n}$ of isotype $\tau$, where
\begin{align}
l  &  =\left\lceil \frac{N-n_{0}+1}{n-1}\right\rceil +1\label{typetau}\\
\tau &  =\left(  n_{0}-1,\left(  n-1\right)  ^{l-2},\tau_{l}\right)
.\nonumber
\end{align}
The number $l=\ell\left(  \tau\right)  $ is the solution of the inequality
$1\leq\tau_{l}=\left(  N-n_{0}+1\right)  -\left(  l-2\right)  \left(
n-1\right)  \leq n-1$ ($\left\lceil r\right\rceil $ denotes the smallest
integer $\geq r$). Then the index for the corresponding singular polynomial is
given by:%
\begin{equation}
\lambda=\left\{
\begin{tabular}
[c]{ll}%
$\left(  m_{0}^{\tau_{2}},0^{n_{0}-1}\right)  ,$ & $l=2$\\
$\left(  \left(  m_{0}+\left(  l-2\right)  m\right)  ^{\tau_{l}},\left(
m_{0}+\left(  l-3\right)  m\right)  ^{n-1},\ldots,m_{0}^{n-1},0^{n_{0}%
-1}\right)  ,$ & $l\geq3$%
\end{tabular}
\right.  . \label{singlb}%
\end{equation}
Note that $l=2$ is equivalent to $N-n_{0}+1<n$ or $d<\frac{n_{0}}{N-n_{0}+1}$,
and $\tau_{2}=N-n_{0}+1$. For $l\geq3$ the computation for $\lambda$ uses the
notation of the previous section with $t_{i}=d+l-1-i,mt_{i}=m_{0}+\left(
l-1-i\right)  m$ for $1\leq i\leq l-1$. The $S_{N}$-module of singular
polynomials is $\mathrm{span}_{\mathbb{Q}}\left\{  w\zeta_{\lambda}^{x}\left(
\kappa_{0}\right)  :w\in S_{N}\right\}  $ and the basis corresponding to
Murphy's construction is exactly the set of $\zeta_{\alpha}^{x}\left(
\kappa_{0}\right)  $ such that $\alpha$ is a reverse lattice permutation of
$\lambda$. There are no other singular polynomials.

The relation of modules of singular polynomials to monodromy representations
of the Hecke algebra was discussed in \cite[Sect.6]{DJO}. The parameter is
$q=e^{-2\pi\mathrm{i}\kappa}$; the existence of singular polynomials of
isotype $\tau$ shows that the monodromy representation corresponding to $\tau$
contains the trivial representation. There is a general result on the
connection between monodromy (called the KZ-functor) and the dual Specht
modules in \cite[Sect. 6.2]{GG}.

Recall the definition of the rational Cherednik algebra (see \cite[Sect.
3]{GG} and \cite{EG}). We consider the image $A\left(  \kappa\right)  $ under
the faithful representation on $\mathcal{P}$; indeed $A\left(  \kappa\right)
$ is the $\mathbb{Q}\left(  \kappa\right)  $-algebra generated by $\left\{
\mathcal{D}_{i}\left(  \kappa\right)  :1\leq i\leq N\right\}  \cup\left\{
x_{i}:1\leq i\leq N\right\}  \cup S_{N}$ (where $x_{i}$ denotes the
multiplication map $p\left(  x\right)  \mapsto x_{i}p\left(  x\right)  $ and
$w\in S_{N}$ acts by $p\left(  x\right)  \mapsto p\left(  xw\right)  $ for
$p\in\mathcal{P}$). In the sequel, when $\kappa$ is specialized to a rational
$\kappa_{0},$we use $\mathcal{P}$ to denote the polynomials with rational
coefficients (that is, $\mathrm{span}_{\mathbb{Q}}\left\{  x^{\alpha}%
:\alpha\in\mathbb{N}_{0}^{N}\right\}  $). Here are some basic results about
$A\left(  \kappa\right)  $-submodules of $\mathcal{P}$.

\begin{proposition}
\label{modsing}Suppose $M$ is a nontrivial proper $A\left(  \kappa_{0}\right)
$-submodule of $\mathcal{P}$ for some $\kappa_{0}\in\mathbb{Q}$, then $M$ is
the direct sum of its homogeneous components $M_{n}:=M\cap\mathcal{P}_{n}$ for
$n\in\mathbb{N}_{0}$, the nonzero component $M_{n_{0}}$ of least degree
($M_{j}=\left\{  0\right\}  $ for $j<n_{0}$) is an $S_{N}$-module of singular
polynomials and $\kappa_{0}$ is a singular value.
\end{proposition}

\begin{proof}
The identity $\sum_{i=1}^{N}x_{i}\mathcal{D}_{i}\left(  \kappa\right)
=\sum_{i=1}^{N}x_{i}\frac{\partial}{\partial x_{i}}+\kappa\sum_{1\leq i<j\leq
N}\left(  1-\left(  i,j\right)  \right)  $ implies that the Euler operator
$\sum_{i=1}^{N}x_{i}\frac{\partial}{\partial x_{i}}\in A\left(  \kappa\right)
$. Hence $M=\sum_{n=0}^{\infty}\left(  M\cap\mathcal{P}_{n}\right)  $. There
exists $n_{0}>0$ such that $M_{n_{0}}\neq\left\{  0\right\}  $ and
$M_{j}=\left\{  0\right\}  $ for $0\leq j<n_{0}$ (or else $M_{0}\neq\left\{
0\right\}  ,1\in M$ and so $M=\mathcal{P}$. Then $\mathcal{D}_{i}\left(
\kappa_{0}\right)  p=0$ for any $p\in M_{n_{0}}$ and $1\leq i\leq N$.
\end{proof}

Say that the \textit{degree }of an $A\left(  \kappa_{0}\right)  $-submodule
$M$ is the least degree of nonzero homogeneous components of $M$, that is,
$\min\left\{  j:M_{j}\neq\left\{  0\right\}  \right\}  $. There is a symmetric
bilinear form on $\mathcal{P}$ defined by
\[
\left\langle p,q\right\rangle _{\kappa}=p\left(  \mathcal{D}_{1}\left(
\kappa\right)  ,\ldots,\mathcal{D}_{N}\left(  \kappa\right)  \right)  q\left(
x\right)  |_{x=0}.
\]
The radical was defined in \cite[Sect. 4]{DJO} to be%
\[
\mathrm{Rad}\left(  \kappa\right)  :=\left\{  p\in\mathcal{P}:\left\langle
p,q\right\rangle _{\kappa}=0\text{ for all }q\in\mathcal{P}\right\}
\]
and was shown to be an $A\left(  \kappa\right)  $-submodule. For $\kappa
_{0}\in\mathbb{Q}$, $\mathrm{Rad}\left(  \kappa_{0}\right)  \neq\left\{
0\right\}  $ exactly when $\kappa_{0}$ is a singular value.

\begin{proposition}
For any singular value $\kappa_{0}$ the radical $\mathrm{Rad}\left(
\kappa_{0}\right)  $ is the largest proper $A\left(  \kappa_{0}\right)
$-submodule of $\mathcal{P}$.
\end{proposition}

\begin{proof}
Suppose $M$ is a nontrivial proper $A\left(  \kappa_{0}\right)  $-submodule.
Suppose $p\in M_{n}=M\cap\mathcal{P}_{n}$ for some $n>0$ and $p\neq0$. Then
for any $q\in\mathcal{P}_{n}$ we have $\left\langle p,q\right\rangle
_{\kappa_{0}}=\left\langle q,p\right\rangle _{\kappa_{0}}=q\left(
\mathcal{D}_{1}\left(  \kappa\right)  ,\ldots,\mathcal{D}_{N}\left(
\kappa\right)  \right)  p\left(  x\right)  \in M_{0}=\left\{  0\right\}  $.
Hence $p\in\mathrm{Rad}\left(  \kappa_{0}\right)  $ and $M\subset
\mathrm{Rad}\left(  \kappa_{0}\right)  $.
\end{proof}

Our complete description of irreducible $S_{N}$-modules of singular
polynomials leads to some explicit results about $A\left(  \kappa_{0}\right)
$-submodules. We use the notation from equations \ref{typetau} and
\ref{singlb}.

\begin{definition}
For any pair $\left(  m_{0},n_{0}\right)  \in\mathbb{N\times}\mathbb{N}$ with
$2\leq n_{0}\leq N$ and $\frac{m_{0}}{n_{0}}\notin\mathbb{N}$ let $M\left(
m_{0},n_{0}\right)  =\left\{  \sum_{i=1}^{n_{\tau}}p_{i}\zeta_{w_{i}\lambda
}^{x}\left(  \kappa_{0}\right)  :p_{i}\in\mathcal{P}\right\}  $, where
$n_{\tau}$ is the degree of the representation $\tau$ and $\left\{
w_{i}\lambda:1\leq i\leq n_{\tau}\right\}  $ is the set of reverse lattice
permutations of $\lambda$ (that is, $\left\{  \zeta_{w_{i}\lambda}^{x}\left(
\kappa_{0}\right)  :1\leq i\leq n_{\tau}\right\}  $ is a basis for the
singular polynomials corresponding to the pair $\left(  m_{0},n_{0}\right)  $).
\end{definition}

The following is from \cite[Sect. 6]{D2}:

\begin{proposition}
$M\left(  m_{0},n_{0}\right)  $ is a proper $A\left(  \kappa_{0}\right)
$-submodule, and its degree is\newline$m\left(  \frac{1}{2}\left(  l-2\right)
\left(  n-1\right)  \left(  2d+l-3\right)  +\tau_{l}\left(  d+l-2\right)
\right)  $ (where $d=\gcd\left(  m_{0},n_{0}\right)  $ $n=n_{0}/d$ and
$m=m_{0}/d$).
\end{proposition}

\begin{proof}
Clearly $M\left(  m_{0},n_{0}\right)  $ is closed under multiplication by
$\mathcal{P}$ and the action of $S_{N}$. Suppose $f=pg$ where $p\in
\mathcal{P}$ and $g\in\mathrm{span}\left\{  \zeta_{w_{i}\lambda}^{x}\left(
\kappa_{0}\right)  :1\leq i\leq n_{\tau}\right\}  $ (that is, $g$ is
singular). By the product rule,%
\[
\mathcal{D}_{i}\left(  \kappa_{0}\right)  f=p\mathcal{D}_{i}\left(  \kappa
_{0}\right)  g+g\left(  \frac{\partial}{\partial x_{i}}p\right)  +\kappa
_{0}\sum_{j\neq i}\left(  \left(  i,j\right)  g\right)  \frac{p\left(
x\right)  -p\left(  x\left(  i,j\right)  \right)  }{x_{i}-x_{j}}%
\]
for $1\leq i\leq N$. But $\mathcal{D}_{i}\left(  \kappa_{0}\right)  g=0$ and
$\frac{p\left(  x\right)  -p\left(  x\left(  i,j\right)  \right)  }%
{x_{i}-x_{j}}$ is a polynomial, thus $\mathcal{D}_{i}\left(  \kappa
_{0}\right)  f\in M\left(  m_{0},n_{0}\right)  $. The degree is $\left\vert
\lambda\right\vert $ (as in equation \ref{singlb}).
\end{proof}

An equivalent formula for $\frac{\left\vert \lambda\right\vert }{m}$ is
$\left(  N-nd+1\right)  \left(  d+l-2\right)  -\frac{1}{2}\left(  n-1\right)
\left(  l-1\right)  \left(  l-2\right)  $. For a given pair $\left(
m,n\right)  $ with $\gcd\left(  m,n\right)  =1$ (and $2\leq n\leq N$) there
are the $A\left(  \kappa_{0}\right)  $-submodules $\left\{  M\left(
dm,dn\right)  :1\leq d\leq\left\lfloor \frac{N}{n}\right\rfloor \right\}  $
($\left\lfloor r\right\rfloor $ denotes the largest integer $\leq r$). The
degree of $M\left(  dm,dn\right)  $ decreases as $d$ increases. This is
obvious because the nonzero part of the index $\lambda$ for $M\left(  \left(
d+1\right)  m,\left(  d+1\right)  n\right)  $ is a substring of the index
$\lambda^{\prime}$ for $M\left(  dm,dn\right)  $; for example take
$N=10,\kappa_{0}=-\frac{1}{3}$ then the values of $\lambda$ (from equation
\ref{singlb}) are $\left(  4^{2},3^{2},2^{2},1^{2},0^{2}\right)  ,~\left(
4,3^{2},2^{2},0^{5}\right)  ,~\left(  3^{2},0^{8}\right)  $ for $d=1,2,3$
respectively. For direct computation, let $l,\tau,\lambda$ and $l^{\prime
},\tau^{\prime},\lambda^{\prime}$ denote the expressions defined in equations
\ref{typetau} and \ref{singlb} for $\left(  m_{0},n_{0}\right)  $ equal to
$\left(  dm,dn\right)  $ and $\left(  \left(  d+1\right)  m,\left(
d+1\right)  n\right)  $ respectively. If $\tau_{l}=1$ then $l^{\prime}=l-2$
and $\tau_{l-2}^{\prime}=n-1$; if $2\leq\tau_{l}\leq n-1$ then $l^{\prime
}=l-1$ and $\tau_{l-1}^{\prime}=\tau_{l}-1$. For both cases $\left\vert
\lambda\right\vert -\left\vert \lambda^{\prime}\right\vert =m\left(
dn+l-2\right)  $. Thus for any given degree of homogeneity there is at most
one irreducible $S_{N}$-module of singular polynomials of that degree (for
$\kappa_{0}=-\frac{m}{n}$). The singular polynomials of least degree
correspond to $\left(  mq,nq\right)  $ where $q=\left\lfloor \frac{N}%
{n}\right\rfloor $.

\begin{proposition}
Suppose $\gcd\left(  m,n\right)  =1$ and $N=nq+r$ with $0\leq r\leq n-1$ (so
$q=\left\lfloor \frac{N}{n}\right\rfloor $) and $k$ denotes the degree of
$\mathrm{Rad}\left(  -\frac{m}{n}\right)  $. If $r<n-1$ then $k=mq\left(
r+1\right)  $, $\mathrm{Rad}\left(  -\frac{m}{n}\right)  \cap\mathcal{P}_{k}$
is of isotype $\left(  nq-1,N-nq+1\right)  ,$ and equals $\mathrm{span}%
_{\mathbb{Q}}\left\{  w\zeta_{\lambda}^{x}\left(  -\frac{m}{n}\right)  :w\in
S_{N}\right\}  $, where $\lambda=\left(  \left(  mq\right)  ^{r+1}%
,0^{nq-1}\right)  $. If $r=n-1$ then $k=m\left(  qn+1\right)  $ and
$\mathrm{Rad}\left(  -\frac{m}{n}\right)  \cap\mathcal{P}_{k}$ is of isotype
$\left(  nq-1,n-1,1\right)  $\ with corresponding $\lambda=\left(  m\left(
q+1\right)  ,\left(  mq\right)  ^{n-1},0^{nq-1}\right)  $.
\end{proposition}

As well as the maximum $A\left(  \kappa_{0}\right)  $-submodule there is a
minimum one, namely, $M\left(  m,n\right)  $.

\begin{proposition}
Suppose $\gcd\left(  m,n\right)  =1$ and $2\leq n\leq N$ then $M\left(
m,n\right)  $ is contained in each nontrivial $A\left(  -\frac{m}{n}\right)  $-submodule.
\end{proposition}

\begin{proof}
Let $M$ be a proper nontrivial $A\left(  -\frac{m}{n}\right)  $-submodule. Let
$s_{0}$ be the degree of $M$. By Proposition \ref{modsing} $M_{s_{0}}$ is an
$S_{N}$-module of singular polynomials, thus $s_{0}=\left\vert \lambda
\right\vert $ for some $\lambda$ given by Equation \ref{singlb}. Suppose
$\lambda$ corresponds to the \ pair $\left(  dm,dn\right)  $ with $1\leq
d\leq\left\lfloor \frac{N}{n}\right\rfloor $, \ then $M\left(  dm,dn\right)
\subset M$. The intersection of any two nontrivial $A\left(  -\frac{m}%
{n}\right)  $-submodules $M_{1}$ and $M_{2}$ is a nontrivial $A\left(
-\frac{m}{n}\right)  $-submodule (if $f\in M_{1}$ and $g\in M_{2}$ then $fg\in
M_{1}\cap M_{2}$). Thus $M\left(  dm,dn\right)  \cap M\left(  m,n\right)  $ is
a nontrivial submodule of $M\left(  m,n\right)  $ which must equal $M\left(
m,n\right)  $ because the degree of $M\left(  m,n\right)  $ is the maximum for
the degrees of $M\left(  dm,dn\right)  $. Hence $M\left(  m,n\right)  \subset
M\left(  dm,dn\right)  \subset M$.
\end{proof}

One could speculate that $\left\{  M\left(  dm,dn\right)  :1\leq
d\leq\left\lfloor \frac{N}{n}\right\rfloor \right\}  $ is the collection of
all nontrivial proper $A\left(  -\frac{m}{n}\right)  $-submodules and that
they are nested, that is, $M\left(  dm,dn\right)  \subset M\left(  \left(
d+1\right)  m,\left(  d+1\right)  n\right)  $. This would be a
characterization of $\mathrm{Rad}\left(  -\frac{m}{n}\right)  $.

\end{document}